\documentclass[preprint,12pt]{elsarticle}

\usepackage{amssymb}
\usepackage{lineno}
\usepackage[dvipsnames,svgnames,x11names]{xcolor}
\usepackage{amsmath}
\numberwithin{equation}{section}
\usepackage[utf8]{inputenc} % allow utf-8 input
\usepackage{hyperref}       % hyperlinks
\usepackage{url}            % simple URL typesetting
\usepackage{booktabs}       % professional-quality tables
\usepackage{amsfonts}       % blackboard math symbols
\usepackage{nicefrac}       % compact symbols for 1/2, etc.
\usepackage{microtype}      % microtypography
\usepackage{cleveref}       % smart cross-referencing
\usepackage{lipsum}         % Can be removed after putting your text content
\usepackage{graphicx}
\usepackage{natbib}
\usepackage{doi}
\usepackage{tikz}
\usepackage{amsthm}
\usepackage{esint}
\tikzset{>=stealth}
\usepackage{bm}
\usepackage{yhmath}
\usetikzlibrary{patterns}
\usepackage{amssymb}
\usepackage{geometry}
\usepackage{threeparttable}
\usepackage{makecell}
\journal{}

\begin{document}

\begin{frontmatter}

%% Title, authors and addresses

%% use the tnoteref command within \title for footnotes;
%% use the tnotetext command for theassociated footnote;
%% use the fnref command within \author or \address for footnotes;
%% use the fntext command for theassociated footnote;
%% use the corref command within \author for corresponding author footnotes;
%% use the cortext command for theassociated footnote;
%% use the ead command for the email address,
%% and the form \ead[url] for the home page:
%% \title{Title\tnoteref{label1}}
%% \tnotetext[label1]{}
%% \author{Name\corref{cor1}\fnref{label2}}
%% \ead{email address}
%% \ead[url]{home page}
%% \fntext[label2]{}
%% \cortext[cor1]{}
%% \affiliation{organization={},
%%             addressline={},
%%             city={},
%%             postcode={},
%%             state={},
%%             country={}}
%% \fntext[label3]{}

\title{Far-field displacement singularity elimination for time-dependent complex variable method on quasi-three dimensional gravitational shallow tunnelling}

%% use optional labels to link authors explicitly to addresses:
%% \author[label1,label2]{}
%% \affiliation[label1]{organization={},
%%             addressline={},
%%             city={},
%%             postcode={},
%%             state={},
%%             country={}}
%%
%% \affiliation[label2]{organization={},
%%             addressline={},
%%             city={},
%%             postcode={},
%%             state={},
%%             country={}}

\author[label1]{Luo-bin Lin}
\ead{luobin_lin@fjut.edu.cn}
\affiliation[label1]{organization={Fujian Provincial Key Laboratory of Advanced Technology and Informatization in Civil Engineering, College of Civil Engineering, Fujian University of Technology},%Department and Organization
  addressline={No. 69 Xueyuan Road, Shangjie University Town}, 
  city={Fuzhou},
  postcode={350118}, 
  state={Fujian},
  country={China}}

\author[label2]{Fu-quan Chen\corref{cor1}}
\cortext[cor1]{Corresponding author}
\ead{phdchen@fzu.edu.cn}
\affiliation[label2]{organization={College of Civil Engineering, Fuzhou University},%Department and Organization
  addressline={No. 2 Xueyuan Road, Shangjie University Town}, 
  city={Fuzhou},
  postcode={350108}, 
  state={Fujian},
  country={China}}

\author[label1]{Chang-jie Zheng}
\ead{zcj@fjut.edu.cn}

\author[label1]{Yi-qun Huang}
\ead{yiqunhuang@fjut.edu.cn}

\begin{abstract}
  %% Text of abstract
  This paper identifies the nonzero resultant and consequent unique displacement singularity of time-dependent complex variable method on quasi-three dimensional shallow tunnelling in visco-elastic and gravitational geomaterial. The quasi-three dimensional problem is equivalently simplified into a plane-strain one using a time-dependent coefficient of convergence confinement method to simulate the progressive release of initial stress field. The unique displacement singularity is thereby eliminated by fixing the far-field ground surface to produce corresponding counter-acting force to equilibriate the nonzero resultant to formalize a strict equilibrium mechanical model. The mixed boundaries of fixed far-field ground surface and nearby free segment form a homogenerous Riemann-Hilbert problem with extra constraints of the virtual traction along tunnel periphery, which is simultaneously solved using an iterative linear system with good numerical stability. The mixed boundary conditions along the ground surface in the whole excavation time span are well satisfied, and detailed comparisons with corresponding finite element solution are conducted. The comparison results are in good agreements, and the proposed solution illustrates high efficiency. More discussions are made on excavation rate, viscosity, and solution convergence. A latent paradox is additionally disclosed for objectivity.
\end{abstract}

%% Graphical abstract
% \begin{graphicalabstract}
%   % \includegraphics{grabs}
%   
% \end{graphicalabstract}

%%Research highlights
% \begin{highlights}
% \item New complex variable method on quasi-three dimensional shallow tunnelling is proposed.
% \item Far-field displacement singularity of gravitational shallow tunnelling is eliminated.
% \item Proposed solution is accurate and very efficient by comparisons with finite element solution.
% \end{highlights}

\begin{keyword}
%% keywords here, in the form: keyword \sep keyword

%% PACS codes here, in the form: \PACS code \sep code

%% MSC codes here, in the form: \MSC code \sep code
%% or \MSC[2008] code \sep code (2000 is the default)
Gravitational shallow tunnelling \sep Complex variable method \sep Far-field displacement singularity \sep Quasi-three dimensional effect \sep Stress and displacement
\end{keyword}

\end{frontmatter}

% \linenumbers

%% main text
\section{Introduction}
\label{sec:introduction}

Gravitational shallow tunnelling is common in geo-engineering. The geomaterial removal within excavation region for gravitational shallow tunnelling would redistribute the stress in the remaining geomaterial, and the original traction along tunnel periphery caused by the initial stress field should be mechanically cancelled by corresponding opposite one, which could be summed up as an upward resultant located within excavation region. Such an upward resultant would brings logarithmic items in Muskhelishvili's complex potentials \cite{Muskhelishvili1966}, which would lead to a unique displacement singularity at infinity in the remaining lower half geomaterial using corresponding complex variable solution. The logarithmic items of complex potentials have been proven to cause a displacement singularity at infinity of the lower half geomaterial in plane strain conditions for circular tunnels \cite{Strack2002phdthesis,Strack_Verruijt2002buoyancy,Verruijt_Strack2008buoyancy,Lu2016}, which clearly violates tunnel engineering facts.

To cancel the unreasonable displacement singularity, several effects have been made. Lin et al. \cite{Self2020JEM} apply a post-processing mathematical strategy by using symmetrical modification and artificially introducing a parameter related to geomaterial depth, however, further uncertainty of the whole solution is laterally caused. To exclude such an uncertainty, Lin et al. \cite{LIN2024appl_math_model} recently put forward a new equilibrium mechanical model to upgrade the original unbalanced mechanical ones \cite{Strack2002phdthesis,Strack_Verruijt2002buoyancy,Verruijt_Strack2008buoyancy,Lu2016} by fixing the far-field ground surface to produce necessary constraining force to equilibriate the unbalanced one distributed along tunnel periphery, and the displacement singularity at infinity has been resolved for circular shallow tunnelling. Furthermore, displacement singularities of noncircular shallow tunnelling considering gravitational gradient \cite{LIN2024comput_geotech_1,lin2024over-under-excavation} are also resolved using similar equilibrium mechanical models.

Though the above-mentioned equilibrium mechanical models have resolved the displacement singularity at infinity of shallow tunnelling in lower half gravitational geomaterial and thereby provide more reasonable stress and displacement in geomaterial, these solutions \cite{LIN2024appl_math_model,LIN2024comput_geotech_1,lin2024over-under-excavation} still remain in elastic plain strain conditions. Meanwhile, real-world tunnelling is three-dimensional, which can be mechanically simulated by time-dependent and quasi-three dimensional models using convergence confinement method \cite{carranza2000application,Unlu2003three_diminsional,vlachopoulos09:_improv,CarranzaTorres2013,paraskevopoulou18:_analy_conver_confin_method}. Several existing time-dependent complex variable solutions \cite{Wang2017deep_twin,Self2020APM_deep_single,Self2021APM_deep_twin,Self2021APM_shallow_twin} have applied relaxation coefficients of convergence confinement method to simulate the quasi-three dimensional effect of tunnel face. Among the above time-dependent complex variable solutions, the one on shallow tunnelling to simulate quasi-three dimensional effect of tunnel face \cite{Self2021APM_shallow_twin} encounters similar displacement singularity at infinity of the elastic plain strain solutions \cite{Strack2002phdthesis,Strack_Verruijt2002buoyancy,Verruijt_Strack2008buoyancy,Lu2016}, indicating a possible neccesity to improve the existing quasi-three dimensional solution.

Therefore, we should wonder the possibility to eliminate the displacement singularity of the time-dependent complex variable solution on quasi-three dimensional shallow tunnelling in gravitational geomaterial by following the equilibirium strategy of fixing the far-field ground. Alternatively, this paper may also serve as a follow-up of our previous study on reasonable mechanical model on shallow tunnelling in gravitational geomaterial \cite{LIN2024appl_math_model} by extending application scenario from the simple and ideal elastic plain strain condition to a more generalized visco-elastic quasi-three dimensional condition.

\section{Problem formation}
\label{sec:problem}

\subsection{Time-dependent mechanical model on shallow tunnelling in gravititational geomaterial}
\label{sec:problem-1}

As shown in Fig. \ref{fig:1}a, a shallow tunnel of depth $ H $ and radius $ R $ is longitudinally excavated in a lower-half gravititational geomaterial of isotropy and small deformation. The shallow tunnelling is time-dependent with a slow excavation rate $ V ({\rm{m/day}}) $ to exclude dynamic stresses, and the three-dimensional effect near tunnel face is approximated by a plain strain one by introducing a dimensionless and time-dependent relaxation coefficient $ U(t) \in (0,1) $ of convergence confinement method \cite{carranza2000application,Unlu2003three_diminsional,vlachopoulos09:_improv,CarranzaTorres2013,paraskevopoulou18:_analy_conver_confin_method} to simulate progressive release of initial stress field near tunnel face (see Fig. \ref{fig:1}b), so that the original three-dimensional excavation problem can be simulated by plain strain ones of continuous time (see Fig. \ref{fig:1}c). Several existing studies have suggested \cite{Wang2017deep_twin,Self2020APM_deep_single,Self2021APM_deep_twin,Self2021APM_shallow_twin} good approximations with measured datum \cite{CarranzaTorres2013,Self2020APM_deep_single} of real three-dimensional simulation.

The time flow to observe excavation should be opposite to the excavation direction. The shallow tunnel is excavated within the time span $ [ t_{0}, \infty ) $, and consists of three sequential stages below. The time span $ [ t_{0}, t_{1} ) $ denotes the phase that the geomaterial is not yet affected by the effect of tunnel face, and is marked as Stage A. The time span $ [ t_{1}, t_{2} ) $ denotes the phase that the preconvergence of tunnel face \cite{Self2020APM_deep_single,Self2021APM_deep_twin,Self2021APM_shallow_twin} starts to affect the geomaterial ahead of tunnel face, and the time moment $ t_{2} $ denotes the time when tunnel face is reached. The time span $ [ t_{2}, t_{3} ) $ denotes the phase that time-dependent convergence occurs due to removal of geomaterial within tunnel periphery, and tunnel face continuously affect the remaining geomaterial \cite{Self2020APM_deep_single,Self2021APM_deep_twin,Self2021APM_shallow_twin} until tunnel convergence is complete. The two phases considering the effect of tunnel face are together marked as Stage B. The time span $ [t_{3}, \infty) $ denotes the phase that the effect of tunnel face is no longer felt, and is marked as Stage C. The mechanical behaviours of these three stages are elaborated below.

(1) In Stage A, no influence of tunnel face is felt, and only the initial stress field is subjected to the to-be-excavated tunnel periphery, as shown in time span $ [ t_{0}, t_{1} ) $. The geomaterial is intact and can be assumed to be elastic with Poisson's ratio $ \nu $ and initial shear modulus $ G_{\infty} + G_{E} $ (see Eq. (\ref{eq:2.7})). The ground surface and to-be-excavated tunnel periphery are denoted $ {\bm{C}}_{1} $ and $ {\bm{C}}_{2} $, respectively. The finite to-be-excavated geomaterial within the contour $ {\bm{C}}_{2} $ is denoted by $ {\bm{D}} $, and $ \overline{\bm{D}} = {\bm{D}} \cup {\bm{C}}_{2} $. The infinite geomaterial beneath contour $ {\bm{C}}_{1} $ is denoted by $ {\bm{\varOmega}}_{0} $, and $ \overline{\bm{\varOmega}}_{0} = {\bm{\varOmega}}_{0} \cup {\bm{C}}_{1} $. The geomaterial between contours $ {\bm{C}}_{1} $ and $ {\bm{C}}_{2} $ is called the remaining geomaterial and denoted by $ {\bm{\varOmega}} $, and $ \overline{\bm{\varOmega}} = {\bm{C}}_{2} \cup {\bm{\varOmega}} \cup {\bm{C}}_{1} $, $ \overline{\bm{\varOmega}}_{0} = \overline{\bm{\varOmega}} \cup {\bm{D}} $. A complex coordinate system $ z = x + {\rm{i}}y $ is established on the ground surface, and the stress field in geomaterial can be expressed as
\begin{equation}
  \label{eq:2.1}
  \left\{
    \begin{aligned}
      & \sigma_{x}^{0}(z) = k_{0} \gamma y \\
      & \sigma_{y}^{0}(z) = \gamma y \\
      & \tau_{xy}^{0}(z) = 0 \\
    \end{aligned}
  \right.
  , \quad 
  \begin{aligned}
    & z = \in \overline{\bm{\varOmega}}_{0} \\
    & t \in [ t_{0}, t_{1} ) \\
  \end{aligned}
\end{equation}
where $ \sigma_{x}^{0}(z) $, $ \sigma_{y}^{0}(z) $, and $ \tau_{xy}^{0}(z) $ denote horizontal, vertical, and shear components of initial stress field, respectively; $ \gamma $ and $ k_{0} $ denote the unit weight and lateral coefficient of the geomaterial, respectively. The traction along to-be-excavated shallow tunnel periphery $ {\bm{C}}_{2} $ can be thereby expressed as
\begin{equation}
  \label{eq:2.2}
  \left\{
    \begin{aligned}
      & X_{i,0}(S) = \sigma_{x}^{0}(S) \cdot \cos \langle \vec{n}_{0},x \rangle + \tau_{xy}^{0}(S) \cdot \cos \langle \vec{n}_{0},y \rangle \\
      & Y_{i,0}(S) = \tau_{xy}^{0}(S) \cdot \cos \langle \vec{n}_{0},x \rangle + \sigma_{y}^{0}(S) \cdot \cos \langle \vec{n}_{0},y \rangle \\
    \end{aligned}
  \right.
  ,\quad S \in {\bm{C}}_{2}
\end{equation}
where $ X_{i,0} $ and $ Y_{i,0} $ denote horizontal and vertical traction components along contour $ {\bm{C}}_{2} $, $ S $ denotes arbitrary point along contour $ {\bm{C}}_{2} $, $ \vec{n}_{0} $ denotes the outward normal vector along the positive induced orientation regarding of region $ \overline{\bm{D}} $ (counter-clockwise), $ \langle \vec{n}_{0},x \rangle $ and $ \langle \vec{n}_{0},y \rangle $ denote the angles between outward normal vector $ \vec{n}_{0} $ and the posititve directions of $ x $ and $ y $ axes, respectively. The path increment along contour $ {\bm{C}}_{2} $ is denoted by $ {\rm{d}}S $, and Eq. (\ref{eq:2.2}) can be rewritten using $ \cos\langle \vec{n}_{0},x \rangle = \frac{{\rm{d}}y}{{\rm{d}}S} $ and $ \cos\langle \vec{n}_{0},y \rangle = -\frac{{\rm{d}}x}{{\rm{d}}S} $ as
\begin{equation}
  \label{eq:2.2'}
  \tag{2.2'}
  \left\{
    \begin{aligned}
      & X_{i,0}(S) = k_{0} \gamma y(S) \frac{{\rm{d}}y(S)}{{\rm{d}}S} \\
      & Y_{i,0}(S) = - \gamma y(S) \frac{{\rm{d}}x(S)}{{\rm{d}}S}
    \end{aligned}
  \right.
  ,\quad S \in {\bm{C}}_{2}
\end{equation}

(2) In Stage B, the effect of tunnel face starts to be felt in the remaining geomaterial $ \overline{\bm{\varOmega}} $, and certain time-dependent additional stress field would be correspondingly generated, whose horizontal, vertical, and shear components are denoted by $ \sigma_{x}^{2}(z,t) $, $ \sigma_{y}^{2}(z,t) $, and $ \tau_{xy}^{2}(z,t) $. The stress in the remaining geomaterial $ \overline{\bm{\varOmega}} $ can subsequently be expressed as
\begin{equation}
  \label{eq:2.3}
  \left\{
    \begin{aligned}
      & \sigma_{x}(z,t) = \sigma_{x}^{0}(z) + \sigma_{x}^{2}(z,t) \\
      & \sigma_{y}(z,t) = \sigma_{y}^{0}(z) + \sigma_{y}^{2}(z,t) \\
      & \tau_{xy}(z,t) = \tau_{xy}^{0}(z) + \tau_{xy}^{2}(z,t) \\
    \end{aligned}
  \right.
  ,\quad
  \begin{aligned}
    & z \in \overline{\bm{\varOmega}} \\
    & t \in [ t_{1}, t_{3} ) \\
  \end{aligned}
\end{equation}
where $ \sigma_{x}(z,t) $, $ \sigma_{y}(z,t) $, and $ \tau_{xy}(z,t) $ denote the horizontal, vertical, and shear components of the stress field in the remaining geomaterial $ \overline{\bm{\varOmega}} $, respectively.

Eq. (\ref{eq:2.3}) indicates that the initial stress field always exists in the remaining geomaterial $ \overline{\bm{\varOmega}} $ in Stage B, so as the initial traction along tunnel periphery (contour $ {\bm{C}}_{2} $) in Eq. (\ref{eq:2.2'}).  Similarly, the additional stress field caused by tunnel face would also generate corresponding additional traction along tunnel periphery (contour $ {\bm{C}}_{2} $) in Fig. \ref{fig:1}c-2, which is called virtual traction \cite{Wang2017deep_twin,Self2020APM_deep_single,Self2021APM_deep_twin,Self2021APM_shallow_twin} and can be expressed as
\begin{equation}
  \label{eq:2.4}
  \left\{
    \begin{aligned}
      & X_{i,2}(S,t) = - U(t) \cdot X_{i,0}(S) = - U(t) \cdot k_{0} \gamma y(S) \frac{{\rm{d}}y(S)}{{\rm{d}}S} \\
      & Y_{i,2}(S,t) = - U(t) \cdot Y_{i,0}(S) = U(t) \cdot \gamma y(S) \frac{{\rm{d}}x(S)}{{\rm{d}}S} \\
    \end{aligned}
  \right.
  ,\quad
  \begin{aligned}
    & S \in {\bm{C}}_{2} \\
    & t \in [ t_{1}, t_{3} ) \\
  \end{aligned}
\end{equation}
where $ X_{i,2}(S,t) $ and $ Y_{i,2}(S,t) $ denote the horizontal and vertical virtual traction components along tunnel periphery, the detailed expression of equivalent coefficient $ U(t) $ can be expressed as \cite{Unlu2003three_diminsional,Self2020APM_deep_single,Self2021APM_deep_twin,Self2021APM_shallow_twin}
\begin{equation}
  \label{eq:2.5}
  U(t) = 
  \left\{
    \begin{aligned}
      & U_{0} + A_{a}\left[ 1 - \exp\left( B_{a} \cdot X(t) \right) \right], \quad t \in [ t_{1}, t_{2} ) \\
      & U_{0} + A_{b}\left[ 1 - \left( \frac{B_{b}}{B_{b}+X(t)} \right)^{2} \right], \quad t \in [ t_{2}, t_{3} ) \\
    \end{aligned}
  \right.
\end{equation}
with
\begin{subequations}
  \label{eq:2.6}
  \begin{equation}
    \label{eq:2.6a}
    \left\{
      \begin{aligned}
        & U_{0} = 0.22\nu + 0.19 \\
        & X(t) = \frac{V}{R}(t - t_{2}) \\
      \end{aligned}
    \right.
  \end{equation}
  \begin{equation}
    \label{eq:2.6b}
    \left\{
      \begin{aligned}
        & A_{a} = - U_{0} \\
        & B_{a} = 0.73\nu + 0.81 \\
      \end{aligned}
    \right.
  \end{equation}
  \begin{equation}
    \label{eq:2.6c}
    \left\{
      \begin{aligned}
        & A_{b} = 1 - U_{0} \\
        & B_{b} = 0.39\nu+0.65 \\
      \end{aligned}
    \right.
  \end{equation}
\end{subequations}
where $ X(t) $ denotes the normalized longitudinal distance to tunnel face. With suitable values of time moments $ t_{1} $ and $ t_{3} $, the values of $ U(t_{1}) \rightarrow 0 $ and $ U(t_{3}) \rightarrow 1 $ would be respectively obtained.

Eq. (\ref{eq:2.6}) explicitly indicates that the Poisson's ratio $ \nu $ is kept constant in Stage B. However, geomaterial near tunnel face would be no more intact, and would be weakened by shallow tunnelling. To exhibit the geomaterial weakening due to excavation, the Poyting-Thomson model for solid material is applied \cite{Self2021APM_shallow_twin} to simulate the decreasing shear modulus of geomaterial as
\begin{equation}
  \label{eq:2.7}
  G(t-t_{1}) =  G_{\infty} + G_{E} \cdot \exp \left[-\frac{G_{E}}{\eta_{E}} \cdot (t - t_{1})\right], \quad t \in [ t_{1}, t_{3} )
\end{equation}
where $ G_{\infty} $ denotes the long-term shear modulus of geomaterial, $ G_{E} $ denotes the spare shear modulus to simulate geomaterial weakening due to tunnel excavation, $ G_{E} + G_{\infty} $ denotes the instantaneous shear modulus before the geomaterial weakening is considered, $ \eta_{E} $ denotes the Maxwell viscosity of the Poyting-Thomson model. Again, we should emphasize that the Poyting-Thomson model in Eq. (\ref{eq:2.3}) only signifies the geomaterial weakening near tunnel due to excavation, which should be uniformly distributed into the whole geomaterial based on the isotropic assumption. Consequently in this paper, the value of $ G_{E} $ should be much smaller than $ G_{\infty} $, and $ -\frac{G_{E}}{\eta_{E}} ({\rm{day}}^{-1}) $ should be also small to denote very slight loss rate of the spare shear modulus uniformly distributed in the geomaterial during Stage B.

(3) In Stage C, the effect of tunnel face is no longer felt in time span $ [ t_{3}, \infty ) $, and the stress in the remaining geomaterial $ \overline{\bm{\varOmega}} $ can be expressed as
\begin{equation}
  \label{eq:2.8}
  \left\{
    \begin{aligned}
      & \sigma_{x}(z,t) = \sigma_{x}^{0}(z) + \sigma_{x}^{\infty}(z) \\
      & \sigma_{y}(z,t) = \sigma_{y}^{0}(z) + \sigma_{y}^{\infty}(z) \\
      & \tau_{xy}(z,t) = \tau_{xy}^{0}(z) + \tau_{xy}^{\infty}(z) \\
    \end{aligned}
  \right.
  ,\quad
  \begin{aligned}
    & z \in \overline{\bm{\varOmega}} \\
    & t \in [ t_{3}, t_{\infty} ) \\
  \end{aligned}
\end{equation}
where $ \sigma_{x}^{\infty}(z) $, $ \sigma_{y}^{\infty}(z) $, and $ \tau_{xy}^{\infty}(z) $ denote horizontal, vertical, and shear components of additional stress field in time span $ [ t_{3}, \infty ) $ in Stage C. The virtual traction along tunnel periphery (contour $ {\bm{C}}_{2} $) in Fig. \ref{fig:1}c can be expressed as
\begin{equation}
  \label{eq:2.9}
  \left\{
    \begin{aligned}
      & X_{i,\infty}(S,t) = - X_{i,0}(S) = - k_{0} \gamma y(S) \frac{{\rm{d}}y(S)}{{\rm{d}}S} \\
      & Y_{i,\infty}(S,t) = - Y_{i,0}(S) = \gamma y(S) \frac{{\rm{d}}x(S)}{{\rm{d}}S} \\
    \end{aligned}
  \right.
  ,\quad
  \begin{aligned}
    & S \in {\bm{C}}_{2} \\
    & t \in [ t_{3}, \infty ) \\
  \end{aligned}
\end{equation}
Eq. (\ref{eq:2.9}) indicates that zero total traction along tunnel periphery is identical to the inner boundary condition in Ref \cite{LIN2024appl_math_model}. Since geomaterial creep would exist for a long time even after the effect of tunnel face, the shear modulus of geomaterial continuously decrease as
\begin{equation}
  \label{eq:2.10}
  G(t-t_{1}) =  G_{\infty} + G_{E} \cdot \exp \left[-\frac{G_{E}}{\eta_{E}} \cdot (t - t_{1})\right], \quad t \in [ t_{3}, \infty )
\end{equation}

\subsection{Stage unification and mechanical model of static equilibrium}%
\label{sub:problem-2}

\begin{figure}[htpb]
  \begin{center}
    \includegraphics[height = \textheight]{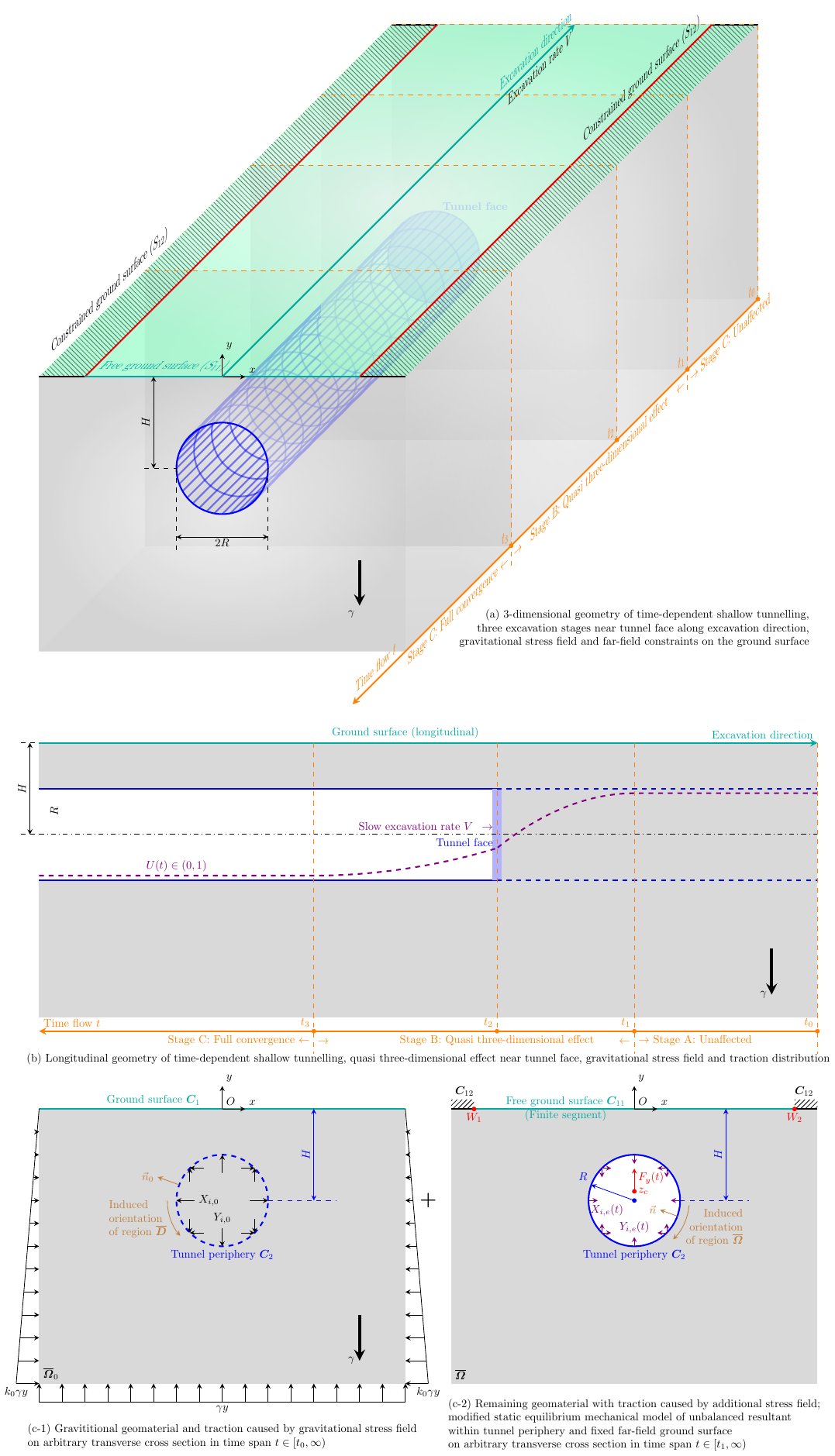}
  \end{center}
  \caption{Schematic diagram of time-dependent shallow tunnelling in gravititional geomaterial}%
  \label{fig:1}
\end{figure}

The equivalent coefficient in Eq. (\ref{eq:2.5}) indicates that $ U(t) \rightarrow 0 $ in time moment $ t_{0} $ in Stage B and $ U(t) \rightarrow 1 $ in time span $ [ t_{3}, \infty ) $ in Stage C, thus, the traction equilibiriums along tunnel periphery (contour $ {\bm{C}}_{2} $) in Fig. \ref{fig:1}c in Eqs. (\ref{eq:2.4}) and (\ref{eq:2.9}) can be unified as
\begin{equation}
  \label{eq:2.11}
  \left\{
    \begin{aligned}
      & X_{i,e}(S,t) = - U(t) \cdot k_{0} \gamma y(S) \frac{{\rm{d}}y(S)}{{\rm{d}}S} \\
      & Y_{i,e}(S,t) = U(t) \cdot \gamma y(S) \frac{{\rm{d}}x(S)}{{\rm{d}}S} \\
    \end{aligned}
  \right.
  ,\quad
  \begin{aligned}
    & S \in {\bm{C}}_{2} \\
    & t \in [ t_{1}, \infty ) \\
  \end{aligned}
\end{equation}
where $ X_{i,e}(S,t) $ and $ Y_{i,e}(S,t) $ denote horizontal and vertical traction components along tunnel periphery (contour $ {\bm{C}}_{2} $) in Fig. \ref{fig:1}c in Stages B and C, respectively. Eq. (\ref{eq:2.11}) precisely substitutes Eq. (\ref{eq:2.4}), while approximately substitutes Eq. (\ref{eq:2.9}). Correspondingly, the stress components in the remaining geomaterial region $ \overline{\bm{\varOmega}} $ in Eqs. (\ref{eq:2.3}) and (\ref{eq:2.8}) can be unified as
\begin{equation}
  \label{eq:2.12}
  \left\{
    \begin{aligned}
      & \sigma_{x}(z,t) = \sigma_{x}^{0}(z) + \sigma_{x}^{e}(z,t) \\
      & \sigma_{y}(z,t) = \sigma_{y}^{0}(z) + \sigma_{y}^{e}(z,t) \\
      & \tau_{xy}(z,t) = \tau_{xy}^{0}(z) + \tau_{xy}^{e}(z,t) \\
    \end{aligned}
  \right.
  ,\quad
  \begin{aligned}
    & z \in \overline{\bm{\varOmega}} \\
    & t \in [ t_{1}, \infty ) \\
  \end{aligned}
\end{equation}
where $ \sigma_{x}^{e}(z,t) $, $ \sigma_{y}^{e}(z,t) $, and $ \tau_{xy}^{e}(z,t) $ denote horizontal, vertical, and shear components of additional stress field in Stages B and C. Eq. (\ref{eq:2.12}) precisely substitutes Eq. (\ref{eq:2.3}), while approximately substitutes Eq. (\ref{eq:2.8}). The shear modulus in Eqs. (\ref{eq:2.7}) and (\ref{eq:2.10}) can also be unified by simple joint of time spans as
\begin{equation}
  \label{eq:2.13}
  G(t-t_{1}) =  G_{\infty} + G_{E} \cdot \exp \left[-\frac{G_{E}}{\eta_{E}} \cdot (t - t_{1})\right], \quad t \in [ t_{1}, \infty )
\end{equation}
Eqs. (\ref{eq:2.11}), (\ref{eq:2.12}), and (\ref{eq:2.13}) unify the boundary conditions, stress distributions, and shear modulus in Stages B and C in mathematical level with sufficient numerical accuracy.

The additional stress components in Eq. (\ref{eq:2.12}) can be expressed via Muskhelishvili's complex potentials \cite{Muskhelishvili1966} as
\begin{equation}
  \label{eq:2.14}
  \left\{
    \begin{aligned}
      & \sigma_{y}^{e}(z,t) + \sigma_{x}^{e}(z,t) = 2 \left[ \frac{\partial\varphi(z,t)}{\partial{z}} + \overline{\frac{\partial\varphi(z,t)}{\partial{z}}} \right] \\
      & \sigma_{y}^{e}(z,t) - \sigma_{x}^{e}(z,t) + 2 {\rm{i}} \tau_{xy}^{e}(z,t) = 2 \left[ \overline{z} \frac{\partial^{2}\varphi(z,t)}{\partial{z}^{2}} + \frac{\partial\psi(z,t)}{\partial{z}} \right] \\
    \end{aligned}
  \right.
  ,\quad
  \begin{aligned}
    & z \in \overline{\bm{\varOmega}} \\
    & t \in [ t_{1}, \infty ) \\
  \end{aligned}
\end{equation}
where $ {\rm{i}} $ denotes imaginary unit with $ {\rm{i}}^{2} = -1 $; $ \varphi(z,t) $ and $ \psi(z,t) $ denote time-dependent complex potentials, which would also cause time-accumulative rectangular displacement in the remaing geomaterial $ \overline{\bm{\varOmega}} $ as \cite{Self2021APM_shallow_twin}
\begin{equation}
  \label{eq:2.15}
  u(z,t) + {\rm{i}}v(z,t) = \frac{1}{2} \int_{t_{1}}^{t} H(t-\tau) \cdot \left[ \kappa \cdot \varphi(z,\tau) - z \frac{\overline{\partial\varphi(z,\tau)}}{\partial{z}} - \overline{\psi(z,\tau)} \right] {\rm{d}}\tau
  , \quad
  \begin{aligned}
    & z \in \overline{\bm{\varOmega}} \\
    & t \in [ t_{1}, \infty ) \\
  \end{aligned}
\end{equation}
where $ u(z,t) $ and $ v(z,t) $ denote horizontal and vertical components of time-accumulative rectangular displacement, respectively; Kolosov parameter $ \kappa = 3 - 4\nu $ for plane strain condition, and
\begin{equation}
  \label{eq:2.16}
  H(t-t_{1}) = \mathcal{L}^{-1} \left[ \frac{1}{T \cdot \hat{G}(T)} \right] = \frac{1}{\eta_{E}} \left( \frac{G_{E}}{G_{E}+G_{\infty}} \right)^{2} \cdot \exp \left[ -\frac{G_{E}}{G_{E}+G_{\infty}} \cdot \frac{G_{\infty}}{\eta_{E}} (t-t_{1}) \right] + \frac{\delta(t-t_{1})}{G_{E}+G_{\infty}}
\end{equation}
wherein $ \mathcal{L}^{-1} $ denotes inverse Laplace transform due to $ \hat{G}(T) = \mathcal{L} [G(t-t_{1})] $, $ \delta(t-t_{1}) $ denotes Dirac delta for time moment $ t_{1} $.

Integrating the traction components in Eq. (\ref{eq:2.11}) along the induced orientation along contour $ {\bm{C}}_{2} $ on the remaing geomaterial $ \overline{\bm{\varOmega}} $ side (clockwise) gives the progressive release of the unbalanced resultant due to shallow tunnelling as
\begin{equation}
  \label{eq:2.17}
  \left\{
    \begin{aligned}
      & F_{x}(t) = \varointclockwise_{{\bm{C}}_{2}} X_{i,e}(S,t) |{\rm{d}}S| = U(t) \cdot k_{0} \gamma \ointctrclockwise_{{\bm{C}}_{2}} y {\rm{d}}y = U(t) \cdot k_{0} \gamma \iint_{\overline{\bm{D}}} 0 \cdot {\rm{d}}x {\rm{d}}y = 0 \\
      & F_{y}(t) = \varointclockwise_{{\bm{C}}_{2}} Y_{i,e}(S,t) |{\rm{d}}S| = - U(t) \cdot \gamma \ointctrclockwise_{{\bm{C}}_{2}} y {\rm{d}}x = U(t) \cdot \gamma \iint_{\overline{\bm{D}}} 1 \cdot {\rm{d}}x {\rm{d}}y = U(t) \cdot \gamma \pi R^{2}
    \end{aligned}
  \right.
\end{equation}
where $ |{\rm{d}}S| = - {\rm{d}}S $ for clockwise integration, and the Green Theorem is applied in the third equalities. Eq. (\ref{eq:2.17}) indicates that a time-dependent nonzero resultant exists along tunnel periphery (contour $ {\bm{C}}_{2} $) in the remaining geomaterial $ \overline{\bm{\varOmega}} $ due to shallow tunnelling in the gravititational geomaterial. This nonzero resultant may be located at arbitrary point $ z_{c} $ within contour $ {\bm{C}}_{2} $, as shown in Fig. \ref{fig:1}c. Therefore, the time-dependent complex potentials in Eqs. (\ref{eq:2.14}) and (\ref{eq:2.15}) containing a unique nonzero resultant can be expanded according to Refs \cite{Strack2002phdthesis,Self2021APM_shallow_twin} as
\begin{equation}
  \label{eq:2.18}
  \left\{
    \begin{aligned}
      & \varphi(z,t) = - U(t) \cdot \frac{{\rm{i}}R^{2}\gamma}{2} \ln(z-\overline{z}_{c}) - U(t) \cdot \frac{{\rm{i}}R^{2}\gamma}{2(1+\kappa)} \ln\frac{z-z_{c}}{z-\overline{z}_{c}} + \varphi_{0}(z,t) \\
      & \psi(z,t) = - U(t) \cdot \frac{{\rm{i}}R^{2}\gamma}{2} \ln(z-\overline{z}_{c}) - U(t) \cdot \frac{{\rm{i}}\kappa R^{2}\gamma}{2(1+\kappa)} \ln\frac{z-z_{c}}{z-\overline{z}_{c}} + \psi_{0}(z,t) \\
    \end{aligned}
  \right.
  ,\quad
  \begin{aligned}
    & z \in \overline{\bm{\varOmega}} \\
    & t \in [ t_{1}, \infty ) \\
  \end{aligned}
\end{equation}
where $ \varphi_{0}(z,t) $ and $ \psi_{0}(z,t) $ denote the analytic and finite components of the time-dependent complex potentials in the remaining geomaterial $ \overline{\bm{\varOmega}} $. The complex potentials in Eq. (\ref{eq:2.18}) implicitly assume that the whole ground surface $ {\bm{C}}_{1} $ is free \cite{Strack2002phdthesis,Self2021APM_shallow_twin}. The complex potentials in Eq. (\ref{eq:2.18}) are derived without considering static equilibrium in the remaining geomaterial, and may be called the "traditional complex potentials". When $ t \rightarrow \infty $, the time-dependent traditional complex potentials would turn to elastic formation, as those in Eq. (2.7) in Ref \cite{LIN2024appl_math_model}.

Eqs. (2.5)-(2.8') in Ref \cite{LIN2024appl_math_model} have shown that the elastic traditional complex potentials containing nonzero resultant would cause displacement singularity at infinity in elastic geomaterial, however, the traditional complex potentials in Eq. (\ref{eq:2.18}) in this paper are time-dependent, and no illustration of such a time-dependent displacement singularity is for sure in a visco-elastic geomaterial yet. Substituting the time-dependent complex potentials in Eq. (\ref{eq:2.18}) into the time-accumulative rectangular displacement in Eq. (\ref{eq:2.15}) yields
\begin{equation}
  \label{eq:2.19}
  \begin{aligned}
    u(z,t) + {\rm{i}}v(z,t) = 
    & \; - \frac{1}{2} \int_{t_{1}}^{t} H(t-\tau) \cdot U(\tau) {\rm{d}}\tau \cdot \frac{{\rm{i}}R^{2}\gamma}{2} \left[ \ln(\overline{z}-\overline{z}_{c}) + \kappa\ln(z-\overline{z}_{c}) \right] \\
    & \; - \frac{1}{2} \int_{t_{1}}^{t} H(t-\tau) \cdot U(\tau) {\rm{d}}\tau \cdot \frac{{\rm{i}}R^{2}\gamma}{2(1+\kappa)} \left( \kappa\ln\frac{z-z_{c}}{z-\overline{z}_{c}} + \ln\frac{\overline{z}-\overline{z}_{c}}{\overline{z}-z_{c}} \right) \\
    & \; - \frac{1}{2} \int_{t_{1}}^{t} H(t-\tau) \cdot U(\tau) {\rm{d}}\tau \cdot \left[ \frac{{\rm{i}}R^{2}\gamma}{2} \frac{z}{\overline{z}-z_{c}}+ \frac{{\rm{i}}R^{2}\gamma}{2(1+\kappa)}\left( \frac{z}{\overline{z}-\overline{z}_{c}} - \frac{z}{\overline{z}-z_{c}} \right) \right] \\
    & \; + \frac{1}{2} \int_{t_{1}}^{t} H(t-\tau) \cdot \left[ \kappa \cdot \varphi_{0}(z,\tau) - z \frac{\overline{\partial\varphi_{0}(z,\tau)}}{\partial{z}} - \overline{\psi_{0}(z,\tau)} \right] {\rm{d}}\tau \\
  \end{aligned}
  , \quad
  \begin{aligned}
    & z \in \overline{\bm{\varOmega}} \\
    & t \in [ t_{1}, \infty ) \\
  \end{aligned}
\end{equation}
Both logarithmic items in the first line of Eq. (\ref{eq:2.19}) would reach a unique singularity when $ z \rightarrow \infty $, indicating that a similar unique displacement singularity at infinity is found in the existing time-dependent complex potentials \cite{Self2021APM_shallow_twin}. Such a displacement singularity obviously violate tunnel engineering fact, and should be eliminated accordingly. In other words, the assumption of free ground surface $ {\bm{C}}_{1} $ for the complex potentials in Eq. (\ref{eq:2.18}) is irrational, and should be modified. 

Similar to the mechanical model in Ref \cite{LIN2024appl_math_model}, we can symmetrically fix the far-field ground surface, as shown in Figs. \ref{fig:1}a and \ref{fig:1}c-2. To be specific, the ground surface $ {\bm{C}}_{1} $ should be separated into a free and finite ground surface segment $ {\bm{C}}_{11} $ above the shallow tunnel and the rest fixed and infinite ground surface segment $ {\bm{C}}_{12} $ in far field, as shown in Fig. \ref{fig:1}c-2. The joint points between segments $ {\bm{C}}_{11} $ and $ {\bm{C}}_{12} $ are denoted by $ W_{1} $ and $ W_{2} $, which should be axisymmetrical about $ y $ axis of rectangular coordinate system $ z = x + {\rm{i}} y $. With the symmetrically fixed far-field ground surface, the vertically time-dependent unbalanced resultant $ F_{y}(t) $ should be always equilibriated by a corresponding equilibrant resultant along the far-field ground surface segment $ {\bm{C}}_{12} $ as
\begin{equation}
  \label{eq:2.20}
  \int_{{\bm{C}}_{12}} \left[ X_{i,e}(W,t) + {\rm{i}} Y_{i,e}(W,t) \right] |{\rm{d}}W| = - {\rm{i}} F_{y}(t) = - U(t) \cdot {\rm{i}} \gamma \pi R^{2}, \quad t \in [ t_{1}, \infty )
\end{equation}
where $ W $ denotes arbitrary point along the fixed far-field ground surface segmenet $ {\bm{C}}_{12} $. The integrating path is from point $ W_{1} $ horizontally towards left to infinity, then from infinity horizontally towards left to point $ W_{2} $, thus, the remaining geomaterial $ \overline{\bm{\varOmega}} $ is always kept on the left side. Such an integrating path is reasonable, since the infinity can be reached from arbitrary path and can reach arbitrary finite point along arbitrary path in complex variable theory. The equilibrium in Eq. (\ref{eq:2.20}) should always exist in the remaining geomaterial in time span $ [ t_{1}, \infty ) $. Simialr equilibriums in elastic geomaterial have been verified in our previous studies \cite{LIN2024appl_math_model,LIN2024comput_geotech_1,lin2024over-under-excavation}, however, the equilibrium in Poyting-Thomson geomaterial in Eq. (\ref{eq:2.20}) is not for sure yet, and would be analytically and numerically examined in Eq. (\ref{eq:4.10c'}) and Section \ref{sec:Case verification}, respectively.

\subsection{Mixed boundary value problem and conformal mapping}%
\label{sub:Mixed boundary value problem}

In time span $ [ t_{1}, \infty ) $, the mechanical model to solve the time-dependent additional stress field in Eq. (\ref{eq:2.14}) and time-accumulative displacement field in Eq. (\ref{eq:2.15}) is the mixed boundary one in Fig. \ref{fig:1}c-2. On the ground surface (unbounded contour $ {\bm{C}}_{1} $), we have
\begin{subequations}
  \label{eq:2.21}
  \begin{equation}
    \label{eq:2.21a}
    X_{i,e}(W,t) + {\rm{i}} Y(W,t) = 0,\quad W \in {\bm{C}}_{11}, t \in [ t_{1}, \infty )
  \end{equation}
  \begin{equation}
    \label{eq:2.21b}
    u(W,t) + {\rm{i}} v(W,t) = 0,\quad W \in {\bm{C}}_{12}, t \in [ t_{1}, \infty )
  \end{equation}
  For the tunnel periphery (bounded contour $ {\bm{C}}_{2} $), the boundary condition in Eq. (\ref{eq:2.11}) can be transformed to its indefinite integral formation for clockwise length increament $ |{\rm{d}}S| = - {\rm{d}}S $ as
  \begin{equation}
    \label{eq:2.21c}
    - {\rm{i}} \int \left[ X_{i,e}(S,t) + {\rm{i}} Y_{i,e}(S,t) \right] |{\rm{d}}S| = U(t) \cdot  \left( - {\rm{i}} k_{0}\gamma \int y(S){\rm{d}}y(S) - \gamma \int y(S){\rm{d}}x(S) \right),\quad S \in {\bm{C}}_{2}, t \in [ t_{1}, \infty )
  \end{equation}
\end{subequations}
The boundary conditions in Eq. (\ref{eq:2.21}) obviously form a mixed boundary value problem.

\begin{figure}[htpb]
  \begin{center}
    \includegraphics[width = \textwidth] {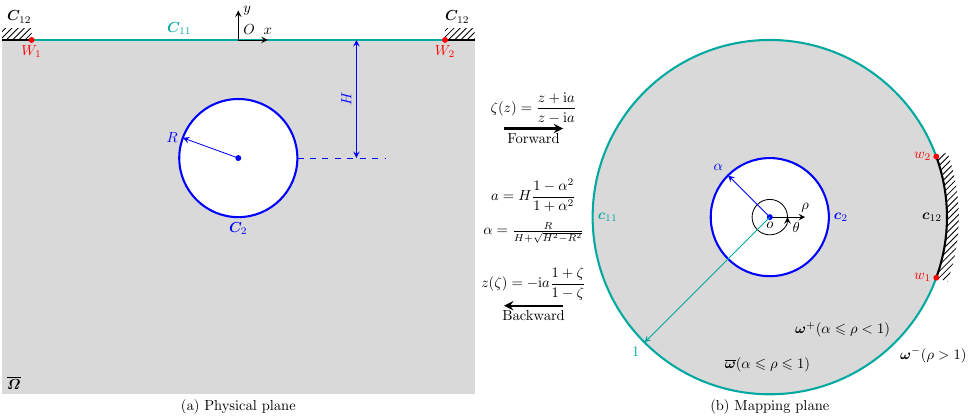}
  \end{center}
  \caption{Bidirectional conformal mapping of the remaining geomaterial}%
  \label{fig:2}
\end{figure}

To solve the above mixed problem via the complex variable method, the Verruijt's conformal mapping \cite{Verruijt1997traction,Verruijt1997displacement} is used to bidirectionally map the remaining geomaterial $ \overline{\bm{\varOmega}} $ onto a unit annulus $ \overline{\bm{\omega}} $ with inner radius $ \alpha $ as
\begin{subequations}
  \label{eq:2.22}
  \begin{equation}
    \label{eq:2.22a}
    \zeta(z) = \frac{z+{\rm{i}}a}{z-{\rm{i}}a},\quad z \in \overline{\bm{\varOmega}}
  \end{equation}
  \begin{equation}
    \label{eq:2.22b}
    z(\zeta) = -{\rm{i}}a\frac{1+\zeta}{1-\zeta}, \quad \zeta \in \overline{\bm{\omega}}
  \end{equation}
  where
  \begin{equation}
    \label{eq:2.22c}
    \left\{
      \begin{aligned}
        & a = H\frac{1-\alpha^{2}}{1+\alpha^{2}} \\
        & \alpha = \frac{R}{H+\sqrt{H^{2}-R^{2}}} \\
      \end{aligned}
    \right.
  \end{equation}
\end{subequations}
The bidirectional conformal mapping is shown in Fig. \ref{fig:2}, where the finite and free ground surface segment $ {\bm{C}}_{11} $, the fixed and infinite ground surface segment $ {\bm{C}}_{12} $, the tunnel periphery $ {\bm{C}}_{2} $, and the joint points $ W_{1} $ and $ W_{2} $ in the remaing geomaterial $ \overline{\bm{\varOmega}} $ are bidirectionally mapped onto the left external arc segment $ {\bm{c}}_{11} $, the right external arc segment $ {\bm{c}}_{12} $, the internal circle $ {\bm{c}}_{2} $, and joint points $ w_{1} $ and $ w_{2} $ in the mapping annulus $ \overline{\bm{\omega}} $, respectively. In addition, $ {\bm{c}}_{1} = {\bm{c}}_{11} \cup {\bm{c}}_{12} $ denotes the whole external periphery of the mapping unit annulus, which is also the mapping of the whole ground surface $ {\bm{C}}_{1} $ in the physical plane. If the joint point coordinates in the physical plane are $ W_{1} = (-x_{0},0) $ and $ W_{2} = (x_{0},0) $, then the corresponding mapping joint point coordinates in the mapping plane would be $ w_{1} = {\rm{e}}^{-{\rm{i}}\theta_{0}} $ and $ w_{2} = {\rm{e}}^{{\rm{i}}\theta_{0}} $, respectively. In other words, the angle range of arc segment $ {\bm{c}}_{12} $ would be $ 2\theta_{0} $, and $ \theta_{0} = -{\rm{i}}\ln\frac{x_{0}+{\rm{i}}a}{x_{0}-{\rm{i}}a} $.

With the bidirectional conformal mapping in Eq. (\ref{eq:2.22}), the curvilinear time-dependent additional stress and displacement field in the remaining geomaterial $ \overline{\bm{\varOmega}} $ in time span $ [ t_{1}, \infty ) $ can be expressed onto the mapping unit annulus $ \overline{\bm{\omega}} $ as
\begin{subequations}
  \label{eq:2.23}
  \begin{equation}
    \label{eq:2.23a}
    \sigma_{\theta}(\zeta,t) + \sigma_{\rho}(\zeta,t) = 2 \left[ \varPhi(\zeta,t) + \overline{\varPhi(\zeta,t)} \right], \quad \zeta \in \overline{\bm{\omega}}, t \in [ t_{1}, \infty )
  \end{equation}
  \begin{equation}
    \label{eq:2.23b}
    \sigma_{\rho}(\zeta,t) + {\rm{i}} \tau_{\rho\theta}(\zeta,t) = \varPhi(\zeta,t) + \overline{\varPhi(\zeta,t)} - \frac{\overline{\zeta}}{\zeta} \left[ \frac{z(\zeta)}{z^{\prime}(\zeta)} \overline{\frac{\partial\varPhi(\zeta,t)}{\partial\zeta}} + \frac{\overline{z^{\prime}(\zeta)}}{z^{\prime}(\zeta)} \overline{\varPsi(\zeta,t)} \right], \quad \zeta \in \overline{\bm{\omega}}, t \in [ t_{1}, \infty )
  \end{equation}
  \begin{equation}
    \label{eq:2.23c}
    g(\zeta,t) = u(\zeta,t) + {\rm{i}} v(\zeta,t) = \frac{1}{2} \int_{t_{1}}^{t} H(t-\tau) \left[ \kappa \cdot \varphi(\zeta,\tau) - z(\zeta)\overline{\varPhi(\zeta,\tau)} - \overline{\psi(\zeta,\tau)} \right]{\rm{d}}\tau, \quad \zeta \in \overline{\bm{\omega}}, t \in [ t_{1}, \infty )
  \end{equation}
  with
  \begin{equation}
    \label{eq:2.23d}
    \left\{
      \begin{aligned}
        & \varPhi(\zeta,t) = \frac{1}{z^{\prime}(\zeta)} \frac{\partial\varphi(\zeta,t)}{\partial\zeta} \\
        & \varPsi(\zeta,t) = \frac{1}{z^{\prime}(\zeta)} \frac{\partial\psi(\zeta,t)}{\partial\zeta} \\
      \end{aligned}
    \right.
  \end{equation}
\end{subequations}
where $ \varphi(\zeta,t) $ and $ \psi(\zeta,t) $ are called the new complex potentials, and should have different formation from and serve as replacement of the traditional ones in Eq. (\ref{eq:2.18}) due to the mixed boundary conditions in Eq. (\ref{eq:2.21}).

The mixed boundary conditions along the ground surface segments $ {\bm{C}}_{11} $ in Eq. (\ref{eq:2.21a}) and $ {\bm{C}}_{12} $ in Eq. (\ref{eq:2.21b}) can be subsequently mapped onto curvilinear ones on the mapping unit annulus $ \overline{\bm{\omega}} $ as
\begin{subequations}
  \label{eq:2.24}
  \begin{equation}
    \label{eq:2.24a}
    \sigma_{\rho}(w,t) + {\rm{i}}\tau_{\rho\theta}(w,t) = 0, \quad w \in {\bm{c}}_{11}, t \in [ t_{1}, \infty )
  \end{equation}
  \begin{equation}
    \label{eq:2.24b}
    u(w,t) + {\rm{i}}v(w,t) = 0, \quad w \in {\bm{c}}_{12}, t \in [ t_{1}, \infty )
  \end{equation}
  where $ w $ denotes arbitrary point along the external periphery $ {\bm{c}}_{1} $ of the mapping unit annulus. Meanwhile, the indefinite intergral formation of traction along tunnel periphery $ {\bm{C}}_{2} $ on the left-hand side of Eq. (\ref{eq:2.21c}) can be rewritten onto its curvilinear formation along the internal circular periphery $ {\bm{c}}_{2} $ of the mapping annulus as
  \begin{equation*}
    \begin{aligned}
      - {\rm{i}} \int \left[ X_{i,e}(S,t) + {\rm{i}} Y_{i,e}(S,t) \right] |{\rm{d}}S| 
      = & \; - {\rm{i}} \int \frac{s}{|s|} \frac{z^{\prime}(s)}{|z^{\prime}(s)|} \left[ \sigma_{\rho}(s,t) + {\rm{i}}\tau_{\rho\theta}(s,t) \right] \cdot |z^{\prime}(s)| \cdot |{\rm{d}}s| \\
      = & \; \int z^{\prime}(s)\left[ \sigma_{\rho}(s,t) + {\rm{i}}\tau_{\rho\theta}(s,t) \right] {\rm{d}}s,\quad s \in {\bm{c}}_{2}, t \in [ t_{1}, \infty )
    \end{aligned}
  \end{equation*}
  where $ |{\rm{d}}s| = - \alpha {\rm{d}}\theta $ for clockwise integration along the internal periphery $ {\bm{c}}_{2} $ in the mapping annulus $ \overline{\bm{\omega}} $. Therefore, the boundary condition along the tunnel periphery $ {\bm{C}}_{2} $ in Eq. (\ref{eq:2.21c}) can be transformed as
  \begin{equation}
    \label{eq:2.24c}
    \int z^{\prime}(s)\left[ \sigma_{\rho}(s,t) + {\rm{i}}\tau_{\rho\theta}(s,t) \right] {\rm{d}}s = U(t) \cdot  \left( - {\rm{i}} k_{0}\gamma \int y(s){\rm{d}}y(s) - \gamma \int y(s){\rm{d}}x(s) \right),\quad s \in {\bm{c}}_{2}, t \in [ t_{1}, \infty )
  \end{equation}
  where the variable of $ x $ and $ y $ on the right-hand side is changed from $ S \in {\bm{C}}_{2} $ in the physical plane to corresponding $ s \in {\bm{c}}_{2} $ in the mapping plane using the backward conformal mapping in Eq. (\ref{eq:2.22b}). 
\end{subequations}

Till now, the mixed boundary condtions in the remaining geomaterial in time span $ [ t_{1}, \infty ) $ in Eq. (\ref{eq:2.21}) are mapped onto corresponding ones in the mapping annulus in Eq. (\ref{eq:2.24}), which can be further expressed by the complex potentials in Eq. (\ref{eq:2.23}) to find solution.

\section{Riemann-Hilbert problem transformation}%
\label{sec:Riemann-Hilbert problem transformation}

The mixed boundary value problem in Eq. (\ref{eq:2.24}) can be transformed to a corresponding Riemann-Hilbert problem to seek solution. To facilitate further discussion, the open regions $ \alpha \leqslant \rho < 1 $ and $ 1 < \rho \leqslant \alpha^{-1} $ in the mapping plane in Fig. \ref{fig:2}b are called $ {\bm{\omega}}^{+} $ and $ {\bm{\omega}}^{-} $, respectively. Clearly, we have $ \overline{\bm{\omega}} = {\bm{\omega}}^{+} \cup {\bm{c}}_{1} $. All the elements on the right-hand sides in Eq. (\ref{eq:2.23}), including new complex potentials and backward conformal mapping, should be analytic within the closure $ \overline{\bm{\omega}} $, according to the mechanical requirement of definition domain $ \zeta \in \overline{\bm{\omega}} $ in Eq. (\ref{eq:2.23}).

Substituting Eq. (\ref{eq:2.23b}) into Eq. (\ref{eq:2.24a}) gives
\begin{equation}
  \label{eq:3.1}
  \varPhi(w,t) = -\overline{\varPhi}(w^{-1},t) + w^{-2} \left[ \frac{z(w)}{z^{\prime}(w)} \frac{\partial\overline{\varPhi}(w^{-1},t)}{\partial w} + \frac{\overline{z}^{\prime}(w^{-1})}{z^{\prime}(w)} \overline{\varPsi}(w^{-1}) \right], \quad w \in {\bm{c}}_{11}, t \in [ t_{1}, \infty )
\end{equation}
where $ \overline{w} = \overline{{\rm{e}}^{{\rm{i}}\theta}} = {\rm{e}}^{-{\rm{i}}\theta} = w^{-1} $ is used; $ \overline{f}(\bullet) $ denotes that only the coefficients of function $ f $ is conjugated, while the independent variable $ \bullet $ is kept as is. If we substitute $ w $ by $ \zeta = \rho {\rm{e}}^{{\rm{i}}\theta} $ with $ 1 < \rho \leqslant \alpha^{-1} $, Eq. (\ref{eq:3.1}) would become
\begin{equation}
  \label{eq:3.2}
  \varPhi(\zeta,t) = -\overline{\varPhi}(\zeta^{-1},t) + \zeta^{-2} \left[ \frac{z(\zeta)}{z^{\prime}(\zeta)} \frac{\partial\overline{\varPhi}(\zeta^{-1},t)}{\partial \zeta} + \frac{\overline{z}^{\prime}(\zeta^{-1})}{z^{\prime}(\zeta)} \overline{\varPsi}(\zeta^{-1},t) \right], \quad \zeta \in {\bm{\omega}}^{-}, t \in [ t_{1}, \infty )
\end{equation}
where $ \zeta^{-1} = \rho^{-1} {\rm{e}}^{-{\rm{i}}\theta} (1 < \rho \leqslant \alpha^{-1}) $ would be within region $ {\bm{\omega}}^{+} $, and $ \overline{\varPhi}(\zeta^{-1},t) $ and $ \overline{\varPsi}(\zeta^{-1},t) $ should be analytic. Therefore, Eq. (\ref{eq:3.2}) extends the definition domain of new complex potential $ \varPhi(\zeta) $ form $ \overline{\bm{\omega}} $ (see Eq. (\ref{eq:2.23b})) to $ \overline{\bm{\omega}} \cup {\bm{\omega}}^{-} $ via analytic continuation.

Eq. (\ref{eq:3.2}) can be modified as
\begin{equation}
  \label{eq:3.3}
  z^{\prime}(\zeta)\varPsi(\zeta,t) = \zeta^{-2}\overline{z}^{\prime}(\zeta^{-1})\overline{\varPhi}(\zeta^{-1},t) + \zeta^{-2}\overline{z}^{\prime}(\zeta^{-1})\varPhi(\zeta,t) - \overline{z}(\zeta^{-1})\frac{\partial\varPhi(\zeta,t)}{\partial\zeta}, \quad \zeta \in \overline{\bm{\omega}}, t \in [ t_{1}, \infty )
\end{equation}
where the complex potentials $ \varPhi(\zeta,t) $ and $ \varPhi(\zeta^{-1},t) $ on the right-hand side are analytic within the regions $ {\bm{\omega}}^{+} $ and $ {\bm{\omega}}^{-} $, respectively, for $ \zeta \in {\bm{\omega}}^{+} $; and both potentials are also analytic along external periphery $ {\bm{c}}_{1} $ of the mapping unit annulus $ \overline{\bm{\omega}} $. According to Eq. (\ref{eq:2.23d}), Eq. (\ref{eq:3.3}) can be rewritten as
\begin{equation}
  \label{eq:3.4}
  \frac{\partial\psi(\zeta,t)}{\partial\zeta} = \partial\left[ -\varphi(\zeta^{-1},t) - \frac{\overline{z}(\zeta^{-1})}{z^{\prime}(\zeta)} \frac{\partial\varphi(\zeta,t)}{\partial\zeta} \right]/{\partial\zeta}, \quad \zeta \in \overline{\bm{\omega}}, t \in [ t_{1}, \infty )
\end{equation}
Conjugate of Eq. (\ref{eq:3.3}) can also be rewritten as
\begin{equation}
  \label{eq:3.5}
  \overline{\varPhi(\zeta,t)} = -\varPhi(\overline{\zeta}^{-1},t) + \overline{\zeta}^{2} \frac{z(\overline{\zeta}^{-1})}{z^{\prime}(\overline{\zeta}^{-1})}\overline{\frac{\partial\varPhi(\zeta,t)}{\partial\zeta}} + \overline{\zeta}^{2}\frac{\overline{z^{\prime}(\zeta)}}{z^{\prime}(\overline{\zeta}^{-1})}\overline{\varPsi(\zeta,t)}, \quad \zeta \in \overline{\bm{\omega}}, t \in [ t_{1}, \infty )
\end{equation}

With substitution of Eq. (\ref{eq:3.5}), Eq. (\ref{eq:2.23b}) can be written into an alternative formation as
\begin{equation}
  \label{eq:3.6}
  \begin{aligned}
    z^{\prime}(\zeta) \left[ \sigma_{\rho}(\zeta,t) + {\rm{i}}\tau_{\rho\theta}(\zeta,t) \right]
    = & \; - \left[ \frac{\overline{\zeta}}{\zeta}z(\zeta) - \overline{\zeta}^{2} z(\overline{\zeta}^{-1}) \frac{z^{\prime}(\zeta)}{z^{\prime}(\overline{\zeta}^{-1})} \right]\overline{\frac{\partial\varPhi(\zeta,t)}{\partial\zeta}} \\
      & \; - \left[ \frac{\overline{\zeta}}{\zeta}\overline{z^{\prime}(\zeta)} - \overline{\zeta}^{2} \overline{z^{\prime}(\zeta)} \frac{z^{\prime}(\zeta)}{z^{\prime}(\overline{\zeta}^{-1})} \right]\overline{\varPsi(\zeta,t)} \\
      & \; + z^{\prime}(\zeta) \varPhi(\zeta,t) - z^{\prime}(\zeta) \varPhi(\overline{\zeta}^{-1},t) 
  \end{aligned}
  , \quad \zeta \in \overline{\bm{\omega}}, t \in [ t_{1}, \infty )
\end{equation}
Taking first derivative of Eq. (\ref{eq:2.23c}) against $ \zeta $ yields
\begin{equation}
  \label{eq:3.7}
  \begin{aligned}
    \frac{\partial g(\zeta,t)}{\partial\zeta} = \frac{1}{2} \int_{t_{1}}^{t} H(t-\tau) \left[ \kappa z^{\prime}(\zeta)\varPhi(\zeta,\tau) - z^{\prime}(\zeta)\overline{\varPhi(\zeta,\tau)} + \frac{\overline{\zeta}}{\zeta}z(\zeta)\overline{\frac{\partial\varPhi(\zeta,\tau)}{\partial\zeta}} + \frac{\overline{\zeta}}{\zeta}\overline{z^{\prime}(\zeta)}\overline{\varPsi(\zeta,\tau)} \right] {\rm{d}}\tau, \quad \zeta \in \overline{\bm{\omega}}, t \in [ t_{1}, \infty ) 
  \end{aligned}
\end{equation}
where $ \frac{{\rm{d}}\overline{\zeta}}{{\rm{d}}\zeta} = \frac{{\rm{d}}\overline{\zeta}}{{\rm{d}}\theta}/\frac{{\rm{d}}\zeta}{{\rm{d}}\theta} = -\frac{\overline{\zeta}}{\zeta} $ with notation $ \zeta = \rho{\rm{e}}^{{\rm{i}}\theta} $ used. With substitution of Eq. (\ref{eq:3.5}), Eq. (\ref{eq:3.7}) can be rewritten into an alternative formation as
\begin{equation}
  \label{eq:3.8}
  \frac{\partial g(\zeta,t)}{\partial\zeta} = \frac{1}{2} \int_{t_{1}}^{t} H(t-\tau) \left\{ 
    \begin{aligned}
      & \; \left[ \frac{\overline{\zeta}}{\zeta}z(\zeta) - \overline{\zeta}^{2} z(\overline{\zeta}^{-1}) \frac{z^{\prime}(\zeta)}{z^{\prime}(\overline{\zeta}^{-1})} \right]\overline{\frac{\partial\varPhi(\zeta,\tau)}{\partial\zeta}} \\
      + & \; \left[ \frac{\overline{\zeta}}{\zeta}\overline{z^{\prime}(\zeta)} - \overline{\zeta}^{2} \overline{z^{\prime}(\zeta)} \frac{z^{\prime}(\zeta)}{z^{\prime}(\overline{\zeta}^{-1})} \right]\overline{\varPsi(\zeta,\tau)} \\
      + & \; \kappa z^{\prime}(\zeta)\varPhi(\zeta,\tau) + z^{\prime}(\zeta) \varPhi(\overline{\zeta}^{-1},\tau)
    \end{aligned}
  \right\} {\rm{d}}\tau, \quad \zeta \in \overline{\bm{\omega}}, t \in [ t_{1}, \infty )
\end{equation}
Eq. (\ref{eq:3.8}) is the first derivative of rectangular displacement field mapped onto the unit annulus $ \overline{\bm{\omega}} $.

When $ \rho \rightarrow 1 $, the items in the brackets $ [] $ in Eqs. (\ref{eq:3.6}) and (\ref{eq:3.8}) would be zero, and these two equations would be simplified as
\begin{subequations}
  \label{eq:3.9}
  \begin{equation}
    \label{eq:3.9a}
    \left. z^{\prime}(\zeta) \left[ \sigma_{\rho}(\zeta,t) + {\rm{i}}\tau_{\rho\theta}(\zeta,t) \right] \right|_{\rho \rightarrow 1} = z^{\prime}(\sigma) \left. \left[ \varPhi(\rho\sigma,t) - \varPhi(\rho^{-1}\sigma,t) \right] \right|_{\rho \rightarrow 1}, \quad t \in [ t_{1}, \infty )
  \end{equation}
  \begin{equation}
    \label{eq:3.9b}
    \left. \frac{\partial g(\zeta,t)}{\partial\zeta} \right|_{\rho \rightarrow 1} = \frac{1}{2} z^{\prime}(\sigma) \int_{t_{1}}^{t} H(t-\tau) \left. \left[ \kappa \varPhi(\rho\sigma,\tau) + \varPhi(\rho^{-1}\sigma,\tau) \right] \right|_{\rho \rightarrow 1} {\rm{d}}\tau, \quad t \in [ t_{1}, \infty )
  \end{equation}
  where $ \sigma = {\rm{e}}^{{\rm{i}}\theta} $. 
\end{subequations}
If we further constrain the polar angle definition domains of Eq. (\ref{eq:3.9}) to approach the curvilinear boundary conditions along the external periphery segments $ {\bm{c}}_{12} $ and $ {\bm{c}}_{11} $ in Eq. (\ref{eq:2.24}), respectively, we have
\begin{subequations}
  \label{eq:3.10}
  \begin{equation}
    \label{eq:3.10a}
    z^{\prime}(\sigma) \left. \left[ \varPhi(\rho\sigma,t) - \varPhi(\rho^{-1}\sigma,t) \right] \right|_{\rho \rightarrow 1} = 0, \quad \sigma \in {\bm{c}}_{11}, t \in [ t_{1}, \infty )
  \end{equation}
  \begin{equation}
    \label{eq:3.10b}
    \frac{1}{2}z^{\prime}(\sigma) \int_{t_{1}}^{t} H(t-\tau) \left. \left[ \kappa \varPhi(\rho\sigma,\tau) + \varPhi(\rho^{-1}\sigma,\tau) \right] \right|_{\rho \rightarrow 1} {\rm{d}}\tau = 0, \quad \sigma \in {\bm{c}}_{12}, t \in [ t_{1}, \infty )
  \end{equation}
  Eq. (\ref{eq:3.10b}) would stand, because it is the first derivative of the curvilinear zero displacement along external periphery segment $ {\bm{c}}_{12} $ of the mapping unit annulus in Eq. (\ref{eq:2.23b}).
\end{subequations}
Eq. (\ref{eq:3.10}) is formally similar to the Riemann-Hilbert problem, except that the time variables in both equations are in different formation. As a matter of fact, the difference of time variable formation in Eq. (\ref{eq:3.10}) can be cancelled by the latent information in Eq. (\ref{eq:2.24c}).

The left-hand side of Eq. (\ref{eq:3.6}) can also be obtained by substitution of Eqs. (\ref{eq:2.23b}) and (\ref{eq:2.23d}) as
\begin{equation}
  \label{eq:3.11}
  z^{\prime}(\zeta) \left[ \sigma_{\rho}(\zeta,t) + {\rm{i}}\tau_{\rho\theta}(\zeta,t) \right] = \partial \left[ \varphi(\zeta,t) + z(\zeta)\overline{\varPhi(\zeta,t)} + \overline{\psi(\zeta,t)} \right] / \partial\zeta, \quad \zeta \in \overline{\bm{\omega}}, t \in [ t_{1}, \infty )
\end{equation}
Substituting Eq. (\ref{eq:3.11}) into the left-hand side of Eq. (\ref{eq:2.24c}) yields
\begin{equation}
  \label{eq:3.12}
  \varphi(s,t) + z(s)\overline{\varPhi(s,t)} + \overline{\psi(s,t)} + C(t) = U(t) \cdot  \left( - {\rm{i}} k_{0}\gamma \int y(s){\rm{d}}y(s) - \gamma \int y(s){\rm{d}}x(s) \right), \quad s \in {\bm{c}}_{2}, t \in [ t_{1}, \infty )
\end{equation}
where $ C(t) $ denotes a time-dependent function. Eq. (\ref{eq:3.12}) shows that the time variable $ t $ only appears in the coefficient $ U(t) $ on right-hand side, which is independent of the space variable $ s $. In other words, the time variable and space variable can be separated. Hence, the left-hand side of Eq. (\ref{eq:3.11}) should also contain the same coefficient $ U(t) $ to ensure such a separation, and the only possibility is to allow the time and space variable separation in the new complex potentials as
\begin{equation}
  \label{eq:3.13}
  \left\{
    \begin{aligned}
      & \varphi(\zeta,t) = U(t) \cdot \varphi(\zeta) \\
      & \psi(\zeta,t) = U(t) \cdot \psi(\zeta) \\
      & C(t) = U(t) \cdot C_{0} \\
    \end{aligned}
  \right.
  , \quad \zeta \in \overline{\bm{\omega}}, t \in [ t_{1}, \infty )
\end{equation}
where $ C_{0} $ is an arbitrary intergral constant.

Substituting Eqs. (\ref{eq:3.13}) and (\ref{eq:2.23d}) into Eqs. (\ref{eq:3.10}) and (\ref{eq:3.12}) yields
\begin{subequations}
  \label{eq:3.14}
  \begin{equation}
    \label{eq:3.14a}
    U(t) \cdot z^{\prime}(\sigma) \left. \left[ \varPhi(\rho\sigma) - \varPhi(\rho^{-1}\sigma) \right] \right|_{\rho \rightarrow 1} = 0, \quad \sigma \in {\bm{c}}_{11}, t \in [ t_{1}, \infty )
  \end{equation}
  \begin{equation}
    \label{eq:3.14b}
    \frac{1}{2}\int_{t_{1}}^{t} H(t-\tau) U(\tau) {\rm{d}}\tau \cdot z^{\prime}(\sigma) \left. \left[ \kappa \varPhi(\rho\sigma) + \varPhi(\rho^{-1}\sigma) \right] \right|_{\rho \rightarrow 1} = 0, \quad \sigma \in {\bm{c}}_{11}, t \in [ t_{1}, \infty )
  \end{equation}
  \begin{equation}
    \label{eq:3.14c}
    U(t) \cdot \left[ \varphi(s) + z(s)\overline{\varPhi(s)} + \overline{\psi(s)} + C_{0} \right] = U(t) \cdot  \left( - {\rm{i}} k_{0}\gamma \int y(s){\rm{d}}y(s) - \gamma \int y(s){\rm{d}}x(s) \right), \quad s \in {\bm{c}}_{2}, t \in [ t_{1}, \infty )
  \end{equation}
\end{subequations}
Since $ U(t) \neq 0 $ according to $ U(t) \in (0,1) $ in Eq. (\ref{eq:2.5}), the time-dependent coefficient and convolution in Eq. (\ref{eq:3.14}) can be cancelled as
\begin{subequations}
  \label{eq:3.15}
  \begin{equation}
    \label{eq:3.15a}
    \varphi^{\prime+}(\sigma) - \varphi^{\prime-}(\sigma) = 0, \quad \sigma \in {\bm{c}}_{11}
  \end{equation}
  \begin{equation}
    \label{eq:3.15b}
    \kappa \varphi^{\prime+}(\sigma) + \varphi^{\prime-}(\sigma) = 0, \quad \sigma \in {\bm{c}}_{11}
  \end{equation}
  \begin{equation}
    \label{eq:3.15c}
    \varphi(s) + \frac{z(s)}{\overline{z^{\prime}(s)}}\overline{\varphi^{\prime}(s)} + \overline{\psi(s)} + C_{0} = - {\rm{i}} k_{0}\gamma \int y(s){\rm{d}}y(s) - \gamma \int y(s){\rm{d}}x(s), \quad s \in {\bm{c}}_{2}
  \end{equation}
\end{subequations}
where $ \varphi^{\prime+}(\sigma) $ and $ \varphi^{\prime-}(\sigma) $ denote the value of $ \varphi^{\prime}(\sigma) = z^{\prime}(\sigma)\varPhi(\sigma) $ approaching from $ {\bm{\omega}}^{+} (\rho < 1) $ and $ {\bm{\omega}}^{-} (\rho > 1) $ sides, respectively. Therefore, Eqs. (\ref{eq:3.15a}) and (\ref{eq:3.15b}) form a homogenerous Riemann-Hilbert problem with extra constraint in Eq. (\ref{eq:3.15c}). Similary, the time-dependent coefficient and convolution in the curvilinear stress and displacement field in the remainig geomaterial $ \overline{\bm{\varOmega}} $ mapped onto the unit annulus in Eq. (\ref{eq:2.23}) can be also separated to obtain corresponding plane strain ones as 
\begin{equation}
  \label{eq:2.23a'}
  \tag{2.23a'}
  \sigma_{\theta}(\zeta) + \sigma_{\rho}(\zeta) = 2 \left[ \varPhi(\zeta) + \overline{\varPhi(\zeta)} \right], \quad \zeta \in \overline{\bm{\omega}}
\end{equation}
\begin{equation}
  \label{eq:2.23b'}
  \tag{2.23b'}
  \sigma_{\rho}(\zeta) + {\rm{i}} \tau_{\rho\theta}(\zeta) = \varPhi(\zeta) + \overline{\varPhi(\zeta)} - \frac{\overline{\zeta}}{\zeta} \left[ \frac{z(\zeta)}{z^{\prime}(\zeta)} \overline{\varPhi^{\prime}(\zeta)} + \frac{\overline{z^{\prime}(\zeta)}}{z^{\prime}(\zeta)} \overline{\varPsi(\zeta)} \right], \quad \zeta \in \overline{\bm{\omega}}
\end{equation}
\begin{equation}
  \label{eq:2.23c'}
  \tag{2.23c'}
  g(\zeta) = \kappa \cdot \varphi(\zeta) - z(\zeta)\overline{\varPhi(\zeta)} - \overline{\psi(\zeta)}, \quad \zeta \in \overline{\bm{\omega}}
\end{equation}
with
\begin{equation}
  \label{eq:2.23d'}
  \tag{2.23d'}
  \left\{
    \begin{aligned}
      & \varPhi(\zeta) = \frac{\varphi^{\prime}(\zeta)}{z^{\prime}(\zeta)} \\
      & \varPsi(\zeta) = \frac{\psi^{\prime}(\zeta)}{z^{\prime}(\zeta)} \\
    \end{aligned}
  \right.
\end{equation}

We may notice that the transformed problem in Eq. (\ref{eq:3.15}) is free from time variable, which is identical to the plane strain reduction by introducing the equivalent coefficient $ U(t) $ in the very begining of Section \ref{sec:problem-1}.

\section{Solution of plane-strain Riemann-Hilbert problem}%
\label{sec:Solution of Riemann-Hilbert problem}

\subsection{Solution preparation}%
\label{sub:Solution preparation}

The general solution of the homogenerous Riemann-Hilbert problem in Eq. (\ref{eq:3.15}) can be expressed via Plemelj formula \cite{Muskhelishvili1966} as
\begin{equation}
  \label{eq:4.1}
  \varphi^{\prime}(\zeta) = X(\zeta) \sum\limits_{n = -\infty}^{\infty} {\rm{i}} f_{n} \zeta^{n}, \quad \zeta \in \overline{\bm{\omega}} \cup {\bm{\omega}}^{-}
\end{equation}
where
\begin{equation}
  \label{eq:4.2}
  X(\zeta) = (\zeta - {\rm{e}}^{-{\rm{i}}\theta_{0}})^{-\frac{1}{2}-{\rm{i}}\lambda} (\zeta - {\rm{e}}^{{\rm{i}}\theta_{0}})^{-\frac{1}{2}+{\rm{i}}\lambda}, \quad \lambda = \frac{\ln\kappa}{2\pi}
\end{equation}
The coefficients $ f_{n} $ should be real owing to axisymmetry about $ y $ axis of the mechanical model in Fig. \ref{fig:1}c-2 \cite{LIN2024appl_math_model}. The function $ X(\zeta) $ can be expanded in regions $ {\bm{\omega}}^{+} $ and $ {\bm{\omega}}^{-} $ via Taylor series as
\begin{subequations}
  \label{eq:4.3}
  \begin{equation}
    \label{eq:4.3a}
    X(\zeta) = \sum\limits_{k=0}^{\infty} \alpha_{k} \zeta^{k}, \quad \zeta \in {\bm{\omega}}^{+}, \quad \alpha_{k} = -{\rm{e}}^{-2\lambda\theta_{0}} \cdot (-1)^{k}\left( c_{k} + \overline{c}_{k} + \sum\limits_{l=1}^{k-1} c_{l}\overline{c}_{k-l} \right)
  \end{equation}
  \begin{equation}
    \label{eq:4.3b}
    X(\zeta) = \sum\limits_{k=1}^{\infty} \beta_{k} \zeta^{-k}, \quad \zeta \in {\bm{\omega}}^{-}, \quad \beta_{k} = (-1)^{k-1} \left( d_{k-1} + \overline{d}_{k-1} + \sum\limits_{l=1}^{k-2} d_{l}\overline{d}_{k-1-l} \right)
  \end{equation}
\end{subequations}
with
\begin{equation*}
  \left\{
    \begin{aligned}
      & c_{0} = \frac{1}{2} \\
      & c_{k} = \prod\limits_{l=1}^{k} \left( \frac{1}{2} - {\rm{i}}\lambda - l \right) \frac{{\rm{e}}^{{\rm{i}}k\theta_{0}}}{k!}, \quad k \geqslant 1 \\
    \end{aligned}
  \right.
  , \qquad
  \left\{
    \begin{aligned}
      & d_{0} = \frac{1}{2} \\
      & d_{k} = \prod\limits_{l=1}^{k} \left( \frac{1}{2} - {\rm{i}}\lambda - l \right) \frac{{\rm{e}}^{-{\rm{i}}k\theta_{0}}}{k!}, \quad k \geqslant 1 \\
    \end{aligned}
  \right.
\end{equation*}

With substitution of Eq. (\ref{eq:4.3}), Eq. (\ref{eq:4.1}) can be expanded as
\begin{subequations}
  \label{eq:4.4}
  \begin{equation}
    \label{eq:4.4a}
    \varphi^{\prime}(\zeta) = \sum\limits_{k=-\infty}^{\infty} {\rm{i}}A_{k}\zeta^{k}, \quad \zeta \in {\bm{\omega}}^{+}, \quad A_{k} = \sum\limits_{n=-k}^{\infty} \alpha_{n+k}f_{-n} 
  \end{equation}
  \begin{equation}
    \label{eq:4.4b}
    \varphi^{\prime}(\zeta) = \sum\limits_{k=-\infty}^{\infty} {\rm{i}}B_{k}\zeta^{k}, \quad \zeta \in {\bm{\omega}}^{-}, \quad B_{k} = \sum\limits_{n=k+1}^{\infty} \beta_{n-k}f_{n} 
  \end{equation}
\end{subequations}
With the variable separation in Eq. (\ref{eq:3.13}), Eq. (\ref{eq:3.4}) can be simplified by substitution of Eq. (\ref{eq:4.4}) as
\begin{equation}
  \label{eq:4.5}
  \psi^{\prime}(\zeta) = - \sum\limits_{k=-\infty}^{\infty} {\rm{i}}B_{-k-2}\zeta^{k} - \left[ \frac{\overline{z}(\zeta^{-1})}{z^{\prime}(\zeta)} \varphi^{\prime}(\zeta) \right]^{\prime}, \quad \zeta \in {\bm{\omega}}^{+}
\end{equation}

Substituting Eq. (\ref{eq:4.4a}) and (\ref{eq:4.5}) into the left-hand side of Eq. (\ref{eq:3.15c}) yields
\begin{equation}
  \label{eq:4.6}
  \begin{aligned}
    \varphi(s) + \frac{z(s)}{\overline{z^{\prime}(s)}}\overline{\varphi^{\prime}(s)} + \overline{\psi(s)} 
    & = {\rm{i}} \sum\limits_{\substack{k=-\infty \\ k \neq 0}}^{\infty} \left( \frac{\alpha^{k}}{k} A_{k-1} - \frac{\alpha^{-k}}{k} B_{k-1} \right) \sigma^{k} + \frac{z(\alpha^{-1}\sigma) - z(\alpha\sigma)}{\overline{z^{\prime}(\alpha\sigma)}} \cdot {\rm{i}} \sum\limits_{k=-\infty}^{\infty} \alpha^{k} A_{k} \sigma^{-k} \\
    & + {\rm{i}}(A_{-1} + B_{-1})\ln\alpha - C_{0} + {\rm{i}}(A_{-1} - B_{-1}) {\rm{Ln}}\sigma
  \end{aligned}
\end{equation}
where $ {\rm{Ln}} $ denotes the multi-valued natural logarithm. Considering the backward conformal mapping in Eq. (\ref{eq:2.22b}), the non-series item above can be expanded as
\begin{equation}
  \label{eq:4.7}
  \frac{z(\alpha^{-1}\sigma) - z(\alpha\sigma)}{\overline{z^{\prime}(\alpha\sigma)}} = \sum\limits_{k=-\infty}^{\infty} e_{k}(\alpha)\sigma^{k}
\end{equation}
where
\begin{equation*}
  \left\{
    \begin{aligned}
      & e_{-k}(\alpha) = 0, \quad k \geqslant 2 \\
      & e_{-1}(\alpha) = -(1-\alpha^{2})\alpha \\
      & e_{k}(\alpha) = (1-\alpha^{2})^{2}\alpha^{k}, \quad k \geqslant 0 \\
    \end{aligned}
  \right.
\end{equation*}
With Eq. (\ref{eq:4.7}), Eq. (\ref{eq:4.6}) can be transformed as
\begin{equation}
  \label{eq:4.6'}
  \tag{4.6'}
  \begin{aligned}
    \varphi(s) + \frac{z(s)}{\overline{z^{\prime}(s)}}\overline{\varphi^{\prime}(s)} + \overline{\psi(s)} 
    & = {\rm{i}} \sum\limits_{\substack{k=-\infty \\ k \neq 0}}^{\infty} \left( \frac{\alpha^{k}}{k} A_{k-1} - \frac{\alpha^{-k}}{k} B_{k-1} + \sum\limits_{l=-\infty}^{\infty} e_{l}(\alpha)\alpha^{l-k} A_{l-k} \right) \sigma^{k} \\
    & + {\rm{i}} \sum\limits_{l=-\infty}^{\infty} e_{l}(\alpha)\alpha^{l}A_{l} + {\rm{i}}(A_{-1} + B_{-1})\ln\alpha - C_{0} + {\rm{i}}(A_{-1} - B_{-1}) {\rm{Ln}}\sigma 
  \end{aligned}
\end{equation}

The components of the right-hand side of Eq. (\ref{eq:3.15c}) can be expanded using the backward conformal mapping in Eq. (\ref{eq:2.22b}) as
\begin{equation*}
  \left\{
    \begin{aligned}
      & x(\alpha\sigma) = -\frac{{\rm{i}}a}{2}\left( \frac{1+\alpha\sigma}{1-\alpha\sigma} - \frac{\sigma+\alpha}{\sigma-\alpha} \right) \\
      & y(\alpha\sigma) = -\frac{a}{2}\left( \frac{1+\alpha\sigma}{1-\alpha\sigma} + \frac{\sigma+\alpha}{\sigma-\alpha} \right) \\
      & \frac{{\rm{d}}x(\alpha\sigma)}{{\rm{d}}\sigma} = -{\rm{i}}a\alpha\left[ \frac{1}{(1-\alpha\sigma)^{2}} + \frac{1}{(\sigma-\alpha)^{2}} \right] \\
      & \frac{{\rm{d}}y(\alpha\sigma)}{{\rm{d}}\sigma} = -a\alpha\left[ \frac{1}{(1-\alpha\sigma)^{2}} - \frac{1}{(\sigma-\alpha)^{2}} \right] \\
    \end{aligned}
  \right.
\end{equation*}
With the above expansions, the integrand of the right-hand side of Eq. (\ref{eq:3.15c}) can be expanded as
\begin{equation}
  \label{eq:4.8}
  -y(\alpha\sigma)\left[ k_{0}\frac{{\rm{d}}y(\alpha\sigma)}{{\rm{d}}\sigma} - {\rm{i}}\frac{{\rm{d}}x(\alpha\sigma)}{{\rm{d}}\sigma} \right] = \sum\limits_{k=-\infty}^{\infty} E_{k} \sigma^{k}
\end{equation}
where the coefficients $ E_{k} $ can be computed using the numerical method of sample point, and the details can be found in Refs \cite{LIN2024comput_geotech_1,lin2024over-under-excavation}. Then the right-hand side of Eq. (\ref{eq:3.15c}) can be obtained by integrating the Laurent series in Eq. (\ref{eq:4.8}) as
\begin{equation}
  \label{eq:4.9}
  - {\rm{i}} k_{0}\gamma \int y(\alpha\sigma){\rm{d}}y(\alpha\sigma) - \gamma \int y(\alpha\sigma){\rm{d}}x(\alpha\sigma) = {\rm{i}} \gamma \int \sum\limits_{k=-\infty}^{\infty} E_{k} \sigma^{k} {\rm{d}}\sigma = {\rm{i}} \sum\limits_{\substack{k=-\infty \\ k \neq 0}}^{\infty} \gamma \frac{E_{k-1}}{k} \sigma^{k} + {\rm{i}} \gamma E_{-1} {\rm{Ln}}\sigma
\end{equation}

According to Eq. (\ref{eq:3.15c}), the coefficients of the same $ \sigma $ order of right-hand sides of Eqs. (\ref{eq:4.6'}) and (\ref{eq:4.9}) should be equal, and we have
\begin{subequations}
  \label{eq:4.10}	
  \begin{equation}
    \label{eq:4.10a}
    A_{-k-1} = \gamma \alpha_{k} E_{-k-1} + \alpha^{2k}B_{-k-1} + k\alpha^{2k}\sum\limits_{l=-\infty}^{\infty} e_{l}(\alpha)\alpha^{l} A_{l+k}, \quad k \geqslant 1
  \end{equation}
  \begin{equation}
    \label{eq:4.10b}
    B_{k-1} = -\gamma\alpha^{k}E_{k-1} + \alpha^{2k}A_{k-1} + k\sum\limits_{l=-\infty}^{\infty} e_{l}(\alpha)\alpha^{l} A_{l-k}, \quad k \geqslant 1
  \end{equation}
  \begin{equation}
    \label{eq:4.10c}
    A_{-1} - B_{-1} = \gamma E_{-1}
  \end{equation}
  The trivial constant items are caused by integration, and are not necessary to present. The value $ \gamma E_{-1} $ in Eq. (\ref{eq:4.10c}) is the imaginary part of the residual of the Laurent series in Eq. (\ref{eq:4.9}), indicating the curvilinear unbalanced resultant along tunnel periphery $ {\bm{c}}_{2} $ on the remaining geomaterial side in the mapping unit annulus $ \overline{\bm{\omega}} $. Furthermore, the displacement single-valuedness requires that \cite{LIN2024appl_math_model,LIN2024comput_geotech_1,lin2024over-under-excavation}
  \begin{equation}
    \label{eq:4.10d}
    \kappa A_{-1} + B_{-1} = 0
  \end{equation}
\end{subequations}

On the other hand, the time-dependent resultant equilibrium in Eq. (\ref{eq:2.20}) should be examined. Since the ground surface segment $ {\bm{C}}_{11} $ is traction-free, we have
\begin{equation}
  \label{eq:2.20'}
  \tag{2.20'}
  \int_{{\bm{C}}_{12}} \left[ X_{i,e}(W,t) + {\rm{i}} Y_{i,e}(W,t) \right] |{\rm{d}}W| = \int_{{\bm{C}}_{1}} \left[ X_{i,e}(W,t) + {\rm{i}} Y_{i,e}(W,t) \right] |{\rm{d}}W|, \quad t \in [ t_{1}, \infty )
\end{equation}
Using the bidirectional conformal mapping, Eq. (\ref{eq:2.20'}) can be transformed as
\begin{equation*}
  \begin{aligned}
    \int_{{\bm{C}}_{1}} \left[ X_{i,e}(W,t) + {\rm{i}} Y_{i,e}(W,t) \right] |{\rm{d}}W| 
    & = \ointctrclockwise_{{\bm{c}}_{1}} \frac{w}{|w|} \frac{z^{\prime}(w)}{|z^{\prime}(w)|} \left[ \sigma_{\rho}(w,t) + {\rm{i}}\tau_{\rho\theta}(w,t) \right] |z^{\prime}(w)| \cdot |{\rm{d}}w|, \quad t \in [ t_{1}, \infty ) \\
    & = - {\rm{i}} \ointctrclockwise_{{\bm{c}}_{1}} z^{\prime}(\sigma) \left[ \sigma_{\rho}(\sigma,t) + {\rm{i}}\tau_{\rho\theta}(\sigma,t) \right] {\rm{d}}\sigma, \quad t \in [ t_{1}, \infty ) \\
  \end{aligned}
\end{equation*}
where $ |{\rm{d}}w| = {\rm{d}}\theta $ for counter-clockwise induced orientation, and $ w = \sigma $. Substituting Eqs. (\ref{eq:3.11}), (\ref{eq:3.13}), (\ref{eq:4.4a}), and (\ref{eq:4.5}) into above equation yields
\begin{equation*}
  - {\rm{i}} \ointctrclockwise_{{\bm{c}}_{1}} z^{\prime}(\sigma) \left[ \sigma_{\rho}(\sigma,t) + {\rm{i}}\tau_{\rho\theta}(\sigma,t) \right] {\rm{d}}\sigma = - {\rm{i}} \ointctrclockwise_{{\bm{c}}_{1}} U(t) \cdot {\rm{i}} \sum\limits_{k=-\infty}^{\infty} \left( A_{k} - B_{-k-2} \right)\sigma^{k} {\rm{d}}\sigma = U(t) \cdot 2\pi {\rm{i}} (A_{-1} - B_{-1})
\end{equation*}
Substituting the above equation back to Eq. (\ref{eq:2.20}) yields
\begin{equation}
  \label{eq:4.10c'}
  \tag{4.10c'}
  A_{-1} - B_{-1} = - \frac{\gamma R^{2}}{2}
\end{equation}
Comparing Eqs. (\ref{eq:4.10c}) and (\ref{eq:4.10c'}), we have $ E_{-1} = - R^{2}/2 $, which should be numerically verified in case verification.

\subsection{Iterative solution}%
\label{sub:Iterative solution}

Eq. (\ref{eq:4.10}) can be used to obtain the solution of $ f_{n} $ in Eq. (\ref{eq:4.1}). Substituting Eqs. (\ref{eq:4.4}) into Eq. (\ref{eq:4.10}), we have
\begin{subequations}
  \label{eq:4.11}	
  \begin{equation}
    \label{eq:4.11a}
    \left\{
      \begin{aligned}
        & \sum\limits_{n=1}^{\infty} \alpha_{n-1}f_{-n} = \frac{\gamma E_{-1}}{1+\kappa} \\
        & \sum\limits_{n=k+1}^{\infty} \alpha_{n-k-1}f_{-n} = \gamma \alpha_{k} E_{-k-1} + \alpha^{2k}B_{-k-1} + k\alpha^{2k}\sum\limits_{l=-\infty}^{\infty} e_{l}(\alpha)\alpha^{l} A_{l+k}, \quad k \geqslant 1 \\
      \end{aligned}
    \right.
  \end{equation}
  \begin{equation}
    \label{eq:4.11b}
    \left\{
      \begin{aligned}
        & \sum\limits_{n=0}^{\infty} \beta_{n+1}f_{n} = -\kappa \frac{\gamma E_{-1}}{1+\kappa} \\
        & \sum\limits_{n=k}^{\infty} \beta_{n-k+1}f_{n} = -\gamma\alpha^{k}E_{k-1} + \alpha^{2k}A_{k-1} + k\sum\limits_{l=-\infty}^{\infty} e_{l}(\alpha)\alpha^{l} A_{l-k}, \quad k \geqslant 1 \\
      \end{aligned}
    \right.
  \end{equation}
\end{subequations}
Eqs. (\ref{eq:4.11a}) and (\ref{eq:4.11b}) respectively constrain $ f_{n} (n \leqslant -1) $ and $ f_{n} (n \geqslant 0) $. 

Assume that $ f_{n} $ can be accumulated as
\begin{equation}
  \label{eq:4.12}
  f_{n} = \sum\limits_{q=0}^{\infty} f_{n}^{(q)}, \quad -\infty < n < \infty
\end{equation}
where $ q $ dentoes iterations. With Eq. (\ref{eq:4.12}), Eq. (\ref{eq:4.11}) can be reorganized in iterative manner \cite{LIN2024appl_math_model,LIN2024comput_geotech_1,lin2024over-under-excavation} as
\begin{subequations}
  \label{eq:4.13}	
  \begin{equation}
    \label{eq:4.13a}
    \left\{
      \begin{aligned}
        & \sum\limits_{n=1}^{\infty} \alpha_{n-1}f_{-n}^{(0)} = \frac{\gamma E_{-1}}{1+\kappa} \\
        & \sum\limits_{n=k+1}^{\infty} \alpha_{n-k-1}f_{-n}^{(0)} = \gamma \alpha^{k} E_{-k-1}, \quad k \geqslant 1 \\
      \end{aligned}
    \right.
  \end{equation}
  \begin{equation}
    \label{eq:4.13b}
    \left\{
      \begin{aligned}
        & \sum\limits_{n=0}^{\infty} \beta_{n+1}f_{n}^{(0)} = -\kappa \frac{\gamma E_{-1}}{1+\kappa} \\
        & \sum\limits_{n=k}^{\infty} \beta_{n-k+1}f_{n}^{(0)} = -\gamma\alpha^{k}E_{k-1}, \quad k \geqslant 1 \\
      \end{aligned}
    \right.
  \end{equation}
\end{subequations}
and
\begin{subequations}
  \label{eq:4.14}
  \begin{equation}
    \label{eq:4.14a}
    \left\{
      \begin{aligned}
        & \sum\limits_{n=1}^{\infty} \alpha_{n-1}f_{-n}^{(q+1)} = 0 \\
        & \sum\limits_{n=k+1}^{\infty} \alpha_{n-k-1}f_{-n}^{(q+1)} = \alpha^{2k}B_{-k-1}^{(q)} + k\alpha^{2k}\sum\limits_{l=-\infty}^{\infty} e_{l}(\alpha)\alpha^{l} A_{l+k}^{(q)}, \quad k \geqslant 1 \\
      \end{aligned}
    \right.
  \end{equation}
  \begin{equation}
    \label{eq:4.14b}
    \left\{
      \begin{aligned}
        & \sum\limits_{n=0}^{\infty} \beta_{n+1}f_{n}^{(q+1)} = 0 \\
        & \sum\limits_{n=k}^{\infty} \beta_{n-k+1}f_{n}^{(q+1)} = \alpha^{2k}A_{k-1}^{(q)} + k\sum\limits_{l=-\infty}^{\infty} e_{l}(\alpha)\alpha^{l} A_{l-k}^{(q)}, \quad k \geqslant 1 \\
      \end{aligned}
    \right.
  \end{equation}
\end{subequations}
where
\begin{equation}
  \label{eq:4.15}
  \left\{
    \begin{aligned}
      & A_{k}^{(q)} = \sum\limits_{n=-k}^{\infty} \alpha_{n+k}f_{-n}^{(q)} \\
      & B_{k}^{(q)} = \sum\limits_{n=k+1}^{\infty} \beta_{n-k}f_{n}^{(q)} \\
    \end{aligned}
  \right.
\end{equation}
The iterations in Eqs. (\ref{eq:4.14}) and (\ref{eq:4.15}) would proceed with $ q : q+1 $. When the threshold $ \max | f_{n}^{(Q+1)} | \leqslant \varepsilon $ is satisfied ($ \varepsilon = 10^{-16} $ for instance), the iteration may stop, and infinite accumulative sum of $ f_{n} $ in Eq. (\ref{eq:4.12}) can be truncated as
\begin{equation}
  \label{eq:4.12'}
  \tag{4.12'}
  f_{n} = \sum\limits_{q=0}^{Q} f_{n}^{(q)}, \quad -\infty < n < \infty
\end{equation}

The iterative solution in Eq. (\ref{eq:4.12})-(\ref{eq:4.12'}) has been proven numerically stable with small condition number \cite{LIN2024appl_math_model}, and has fast convergence within two seconds \cite{LIN2024appl_math_model,LIN2024comput_geotech_1,lin2024over-under-excavation}.

\section{Time-dependent stress and displacement fields}%
\label{sec:Time-dependent stress and displacement fields}

\subsection{Analytical solution}%
\label{sub:Analytical solution}

The solution in Eq. (\ref{eq:4.12'}) would give the complex potentials in Eqs. (\ref{eq:4.4a}) and subsequent (\ref{eq:4.5}), and the curvilinear stress and displacement mapped onto the unit annulus $ \overline{\bm{\omega}} $ in Eqs. (\ref{eq:2.23a'}), (\ref{eq:2.23b'}), and (\ref{eq:2.23c'}) can be respectively obtained \cite{LIN2024comput_geotech_1,lin2024over-under-excavation} as
\begin{subequations}
  \label{eq:5.1}
  \begin{equation}
    \label{eq:5.1a}
    \sigma_{\theta}(\rho\sigma) + \sigma_{\rho}(\rho\sigma) = 4\Re\left[ \frac{{\rm{i}}}{z^{\prime}(\rho\sigma)} \sum\limits_{k=-\infty}^{\infty} A_{k}\rho^{k}\sigma^{k} \right], \quad \alpha \leqslant \rho \leqslant 1
  \end{equation}
  \begin{equation}
    \label{eq:5.1b}
    \sigma_{\rho}(\rho\sigma) + {\rm{i}}\tau_{\rho\theta}(\rho\sigma) = \frac{{\rm{i}}}{z^{\prime}(\rho\sigma)} \sum\limits_{k=-\infty}^{\infty} \left[ A_{k}\rho^{k} - B_{k}\rho^{-k-2} + (k+1)\sum\limits_{l=-\infty}^{\infty} e_{l}(\rho)A_{l-k-1}\rho^{l-k-2} \right] \sigma^{k}, \quad \alpha \leqslant \rho \leqslant 1
  \end{equation}
  \begin{equation}
    \label{eq:5.1c}
    \begin{aligned}
      g(\rho\sigma) 
      = & {\rm{i}} \sum\limits_{\substack{k=-\infty \\ k \neq 0}}^{\infty} \left( \kappa A_{k-1}\rho^{k} + B_{k-1}\rho^{-k} \right) \frac{\sigma^{k}}{k} - {\rm{i}} \sum\limits_{\substack{k=-\infty \\ k \neq 0}}^{\infty} \left( \kappa A_{k-1} + B_{k-1} \right)\frac{1^{k}}{k} \\
        & - {\rm{i}} \sum\limits_{k=-\infty}^{\infty} \sum\limits_{l=-\infty}^{\infty} e_{l}(\rho)A_{l-k}\rho^{l-k}\sigma^{k} + {\rm{i}} (\kappa A_{-1} - B_{-1})\ln\rho, \quad \alpha \leqslant \rho \leqslant 1
    \end{aligned}
  \end{equation}
\end{subequations}
where $ e_{l}(\rho) $ is similarly computed using the expansion in Eq. (\ref{eq:4.7}) by simply replacing $ \alpha $ by $ \rho $. 

The plane strain results in Eq. (\ref{eq:5.1}) should be restored to time-dependent ones by multiplying time-dependent coefficient and convolution according to the separation in Eq. (\ref{eq:3.13}) as
\begin{subequations}
  \label{eq:5.2}
  \begin{equation}
    \label{eq:5.2a}
    \sigma_{\theta}(\zeta,t) + \sigma_{\rho}(\zeta,t) = U(t) \cdot \left[ \sigma_{\theta}(\zeta) + \sigma_{\rho}(\zeta) \right], \quad \zeta \in \overline{\bm{\omega}}, t \in [ t_{1},\infty )
  \end{equation}
  \begin{equation}
    \label{eq:5.2b}
    \sigma_{\rho}(\zeta,t) + {\rm{i}}\tau_{\rho\theta}(\zeta,t) = U(t) \cdot \left[ \sigma_{\rho}(\zeta) + {\rm{i}}\tau_{\rho\theta}(\zeta) \right], \quad \zeta \in \overline{\bm{\omega}}, t \in [ t_{1},\infty )
  \end{equation}
  \begin{equation}
    \label{eq:5.2c}
    u(\zeta,t) + {\rm{i}}v(\zeta,t) = \frac{1}{2} \int_{t_{1}}^{t} H(t-\tau)U(\tau){\rm{d}}\tau \cdot g(\zeta), \quad \zeta \in \overline{\bm{\omega}}, t \in [ t_{1},\infty )
  \end{equation}
\end{subequations}
where
\begin{equation}
  \label{eq:5.3}
  \int_{t_{1}}^{t} H(t-\tau)U(\tau){\rm{d}}\tau = \frac{1}{\eta_{E}} \left( \frac{G_{E}}{G_{E}+G_{\infty}} \right)^{2} \cdot \int_{t_{1}}^{t} \exp \left[ -\frac{G_{E}}{G_{E}+G_{\infty}} \cdot \frac{G_{\infty}}{\eta_{E}} (t-\tau) \right] U(\tau) {\rm{d}}\tau + \frac{U(t)}{G_{E}+G_{\infty}}
\end{equation}
If geomaterial near tunnel is assumed to be intact before and after tunnel excavation, the shear modulus can be reduced to be elastic in the whole shallow tunnelling procedure, namely $ G(t) = G_{\infty} $ by setting the spare shear modulus $ G_{E} = 0 $ in Eq. (\ref{eq:5.3}), then displacement without time accumulation would be obtained as an alternative of Eq. (\ref{eq:5.2c}):
\begin{equation}
  \label{eq:5.2c'}
  \tag{5.2c'}
  u(\zeta,t) + {\rm{i}}v(\zeta,t) = \frac{U(t)}{2G_{\infty}} \cdot g(\zeta), \quad \zeta \in \overline{\bm{\omega}}, t \in [ t_{1},\infty )
\end{equation}

With the curvilinear time-dependent stress and displacement fields mapped onto the unti annulus $ \overline{\bm{\omega}} $, the rectangular stress and displacement fields caused by shallow tunnelling in the remaining geomaterial $ \overline{\bm{\varOmega}} $ can be computed via the bidirectional conformal mapping in Eq. (\ref{eq:2.22b}) as
\begin{subequations}
  \label{eq:5.4}
  \begin{equation}
    \label{eq:5.4a}
    \sigma_{y}^{e}(z,t) + \sigma_{x}^{e}(z,t) = \sigma_{\theta}[\zeta(z),t] + \sigma_{\rho}[\zeta(z),t], \quad z \in \overline{\bm{\varOmega}}, t \in [ t_{1},\infty )
  \end{equation}
  \begin{equation}
    \label{eq:5.4b}
    \sigma_{y}^{e}(z,t) - \sigma_{x}^{e}(z,t) + 2{\rm{i}}\tau_{xy}^{e}(z,t) = \left\{ \sigma_{\theta}[\zeta(z),t] - \sigma_{\rho}[\zeta(z),t] + 2{\rm{i}}\tau_{\rho\theta}[\zeta(z),t] \right\} \cdot \frac{\overline{\zeta(z)}}{\zeta(z)} \left. \frac{\overline{z^{\prime}(\zeta)}}{z^{\prime}(\zeta)} \right|_{\zeta \rightarrow \zeta(z)}, \quad z \in \overline{\bm{\varOmega}}, t \in [ t_{1},\infty )
  \end{equation}
  \begin{equation}
    \label{eq:5.4c}
    u(z,t) + {\rm{i}}v(z,t) = \left. \left[ u(\zeta,t) + {\rm{i}}v(\zeta,t) \right] \right|_{\zeta \rightarrow \zeta(z)}, \quad z \in \overline{\bm{\varOmega}}, t \in [ t_{1},\infty )
  \end{equation}
\end{subequations}
The total rectangular stress field in the remaining geomaterial $ \overline{\bm{\varOmega}} $ can be obtained by Eq. (\ref{eq:2.12}).

\subsection{Numerical results}%
\label{sub:Numerical results}

The time-dependent stress and displacement fields above is analytical and theoretical by setting infinite items of Eq. (\ref{eq:4.1}) to satisfy the traction and displacement jump values along the external periphery $ {\bm{c}}_{1} $ in the homogenerous Riemann-Hilbert problem in Eq. (\ref{eq:3.15}). In numerical and practical computation, Eq. (\ref{eq:4.1}) has to be truncated into $ 2N+1 $ items ($ -N \leqslant n \leqslant N $). Correspondingly, Eq. (\ref{eq:4.4}) would be truncated as
\begin{equation}
  \label{eq:4.4a'}
  \tag{4.4a'}
  \varphi^{\prime}(\zeta) = \sum\limits_{k=-N}^{N} {\rm{i}}A_{k}\zeta^{k}, \quad \zeta \in {\bm{\omega}}^{+}, \quad A_{k} = \sum\limits_{n=-k}^{N} \alpha_{n+k}f_{-n} 
\end{equation}
\begin{equation}
  \label{eq:4.4b'}
  \tag{4.4b'}
  \varphi^{\prime}(\zeta) = \sum\limits_{k=-N}^{N} {\rm{i}}B_{k}\zeta^{k}, \quad \zeta \in {\bm{\omega}}^{-}, \quad B_{k} = \sum\limits_{n=k+1}^{N} \beta_{n-k}f_{n} 
\end{equation}
The results in Eq. (\ref{eq:5.1}) should also be truncated as
\begin{equation}
  \label{eq:5.1a'}
  \tag{5.1a'}
  \sigma_{\theta}(\rho\sigma) + \sigma_{\rho}(\rho\sigma) = 4\Re\left[ \frac{{\rm{i}}}{z^{\prime}(\rho\sigma)} \sum\limits_{k=-N}^{N} A_{k}\rho^{k}\sigma^{k} \right], \quad \alpha \leqslant \rho \leqslant 1
\end{equation}
\begin{equation}
  \label{eq:5.1b'}
  \tag{5.1b'}
  \sigma_{\rho}(\rho\sigma) + {\rm{i}}\tau_{\rho\theta}(\rho\sigma) = \frac{{\rm{i}}}{z^{\prime}(\rho\sigma)} \sum\limits_{k=-N}^{N} \left[ A_{k}\rho^{k} - B_{k}\rho^{-k-2} + (k+1)\sum\limits_{l=-N+k+1}^{N+k+1} e_{l}(\rho)A_{l-k-1}\rho^{l-k-2} \right] \sigma^{k}, \quad \alpha \leqslant \rho \leqslant 1
\end{equation}
\begin{equation}
  \label{eq:5.1c'}
  \tag{5.1c'}
  \begin{aligned}
    g(\rho\sigma) 
    = & {\rm{i}} \sum\limits_{\substack{k=-N+1 \\ k \neq 0}}^{N+1} \left( \kappa A_{k-1}\rho^{k} + B_{k-1}\rho^{-k} \right) \frac{\sigma^{k}}{k} - {\rm{i}} \sum\limits_{\substack{k=-N+1 \\ k \neq 0}}^{N+1} \left( \kappa A_{k-1} + B_{k-1} \right)\frac{1^{k}}{k} \\
      & - {\rm{i}} \sum\limits_{k=-N}^{N} \sum\limits_{l=-N+k}^{N+k} e_{l}(\rho)A_{l-k}\rho^{l-k}\sigma^{k} + {\rm{i}} (\kappa A_{-1} - B_{-1})\ln\rho, \quad \alpha \leqslant \rho \leqslant 1
  \end{aligned}
\end{equation}
Simultaneously, Eqs. (\ref{eq:4.7}) and (\ref{eq:4.8}) should be truncated into $ 2M+1 $ items ( $ -M \leqslant k \leqslant M $ ), and $ M $ should be much larger than $ N $ to ensure numerical accuracy.

Apparently, the finite items of the stress and displacement results above can not fully simulate those jump values in Eq. (\ref{eq:3.15}), and would contain certain numerical errors. Such numerical errors would cause an oscillation pattern of the stress and displacement fields \cite{LIN2024appl_math_model,LIN2024comput_geotech_1,lin2024over-under-excavation}, which can be greatly reduced by applying the Lanczos filtering \cite{LIN2024appl_math_model,LIN2024comput_geotech_1,lin2024over-under-excavation} to replace $ \sigma^{k} $ and $ 1^{k} $ with $ L_{k} \cdot \sigma^{k} $ and $ L_{k} \cdot 1^{k} $ in Eqs. (\ref{eq:5.1a'}), (\ref{eq:5.1b'}), and (\ref{eq:5.1c'}), respectively, where
\begin{equation}
  \label{eq:5.5}
  L_{k} = 
  \left\{
    \begin{aligned}
      & 1, \quad k = 0 \\
      & \sin\left( \frac{k}{N}\pi \right)/\left( \frac{k}{N}\pi \right), \quad {\rm{otherwise}}
    \end{aligned}
  \right.
\end{equation}
Subsequent computation in Eqs. (\ref{eq:5.2}) and (\ref{eq:5.3}) can proceed to obtain numerical results. The convolution in Eq. (\ref{eq:5.3}) can be approximately computed using numerical integration for simplicity.

\section{Case verification}%
\label{sec:Case verification}

In this section, a numerical case is conducted to verify the proposed analytical solution by examining the mixed boundary conditions along ground surface and comparing with a corresponding finite element solution. To be consistent with the plain strain base of the time-dependent analytical solution, the corresponding finite element solution should be conducted using two-dimensional model with plain strain elements. The analytical solution is conducted using programming code \texttt{Fortran} of compiler \texttt{GCC-14.1.1}, and the real and complex linear systems in the solution are solved using \texttt{dgesv} and \texttt{zgesv} packages, respetively. All the data figures are plotted using \texttt{Gnuplot-6.0}.

\subsection{Parameter selection}%
\label{sub:Parameter selection}

The parameters for case verification are listed in Table \ref{tab:1}. According to parameter assumptions in Section \ref{sec:problem-1}, the spare shear modulus $ G_{E} $ is much smaller than the long-term one $ G_{\infty} $ and the viscosity $ \eta_{E} $ takes a large value to simulate the very slight weakening of geomaterial due to excavation. Refs \cite{LIN2024appl_math_model, LIN2024comput_geotech_1, lin2024over-under-excavation} suggest a large $ x_{0} $ to guarantee convergence and mechanical improvement of the present solution, however, the value $ x_{0} $ here takes an irrationally small value for a particular reason to highlight the acute change of both stress and displacement in the computing results for better illustration and comparison. Excavation rate takes $ V = 2{\rm{m/day}} $ to simulate a relatively slow shallow tunnelling for good observation of the effect of tunnel face. The time span $ [t_{0}, t_{1}) $ takes 100 days to simulate sufficient time for initial stress field before the excavation is felt. The time span $ [ t_{1}, t_{2} ) $ and $ [t_{2}, t_{4}] $ would provide sufficient time spans to reveal the effect of tunnel face. The truncation numbers $ N = 200 $ and $ M = 500 $ are selected for good solution convergence. 

\begin{table}[htpb]
  \centering
  \caption{Parameters for case verification}
  \label{tab:1}
  \small
  \begin{tabular}{ccccccccc}
    \toprule
    $ k_{0} $ & $ \gamma ({\rm{kN/m^{3}}}) $ & $ \nu $ & $ G_{\infty} ({\rm{MPa}}) $ & $ G_{E} ({\rm{MPa}}) $ & $ \eta_{E} ({\rm{MPa \cdot day}}) $ & $ x_{0} ({\rm{m}}) $ & $ R ({\rm{m}}) $ & $ H ({\rm{m}}) $ \\
    \midrule
    0.8 & 20 & 0.3 & 20 & 1 & 100000 & 10 & 5 & 10 \\
    \bottomrule
		$ V ({\rm{m/day}}) $ & $ t_{0} ({\rm{day}}) $ & $ t_{1} ({\rm{day}}) $ & $ t_{2} ({\rm{day}}) $ & $ t_{3} ({\rm{day}}) $ & $ t_{4} ({\rm{day}}) $ &  & $ N $ & $ M $ \\
    \midrule
		2 & 0 & 100 & 105 & 110 & 120 &  & 200 & 500 \\
    \bottomrule
  \end{tabular}
\end{table}

With the parameters in Table \ref{tab:1}, the equivalent coefficient of tunnel face in Eq. (\ref{eq:2.5}) can be specified as
\begin{equation}
  \label{eq:2.5'}
  \tag{2.5'}
  U(t) = 
  \left\{
    \begin{aligned}
      & 0.256 \cdot \exp \left[ 0.4116 \cdot \left( t - 105 \right) \right], \quad t \in [ 100, 105 ) \\
      & 1 - 0.744 \cdot \left[ \frac{0.767}{0.767+0.4 \cdot (t-105)} \right]^{2}, \quad t \in [105, 120] \\
    \end{aligned}
  \right.
\end{equation}
From Eq. (\ref{eq:2.5'}), we have $ U(t_{1}) \approx 0.0327 $, $ U(t_{2}) = 0.256 $, $ U(t_{3}) \approx 0.9427 $, and $ U(t_{4}) \approx 0.9923 $, indicating the effect of tunnel face is not felt before and after time moments $ t_{1} $ and $ t_{3} $, respectively. The incremental time variable of numerical integration to compute the convolution in Eq. (\ref{eq:5.3}) takes $ {\rm{d}}\tau \rightarrow 0.01 $ to meet accuracy requirement, and the excavation time span $ [t_{1},t_{4}] $ would contain $ \frac{t_{4}-t_{1}}{{\rm{d}}\tau} = 20000 $ time moments.

\subsection{Examination of mixed boundary conditions}%
\label{sub:Examination of mixed boundary conditions}

\begin{figure}[htpb]
  \centering
  \includegraphics[width = 1.\textwidth]{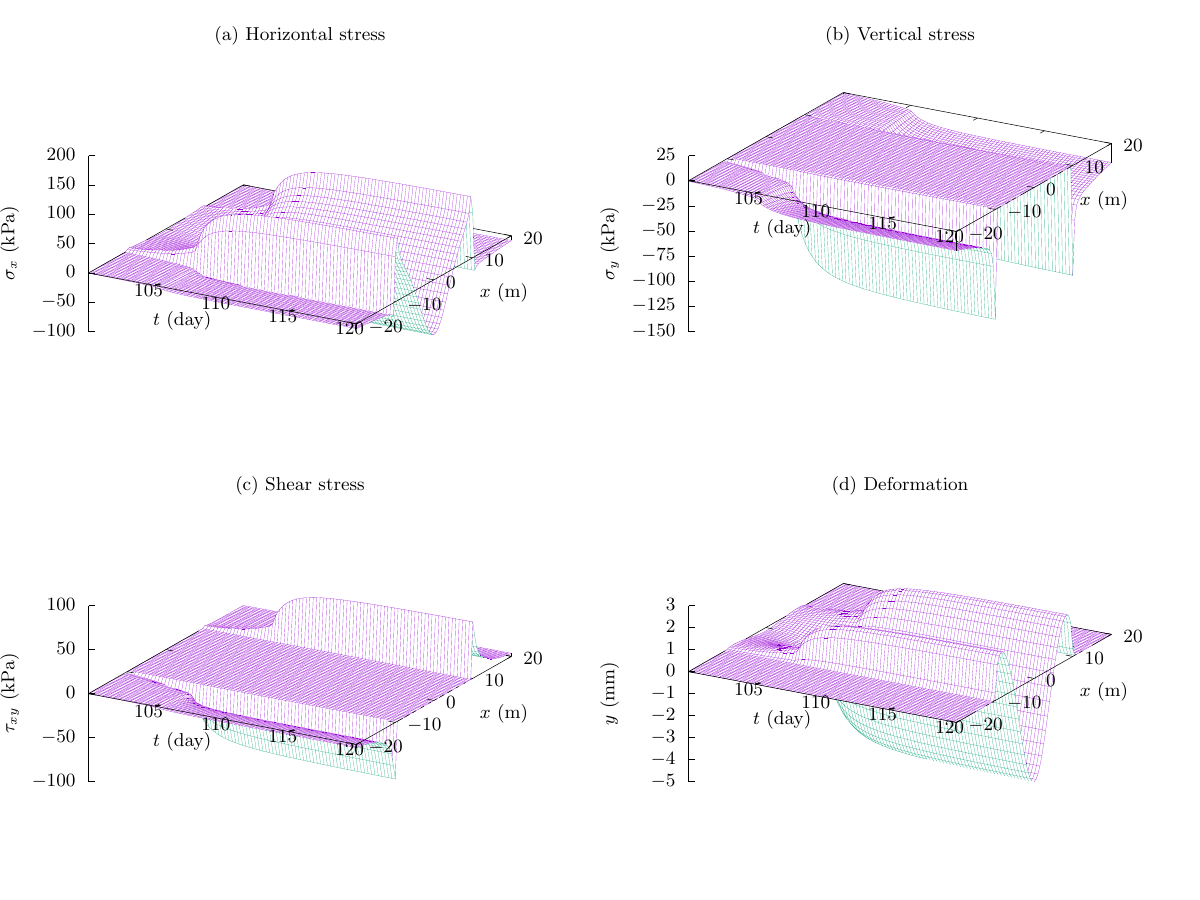}
  \caption{Stress components and deformation on the ground surface against time}
  \label{fig:3}
\end{figure}

The mixed boundary conditions in Eq. (\ref{eq:2.21}) would be examined first for its acute change along ground surface and clear zero values for the whole time span $ [t_{1}, t_{4}] $. Substituting the parameters into the present solution visually gives the stress components and deformation on the ground surface against the whole time span $ [t_{1}, t_{4}] $ in Fig. \ref{fig:3}. Figs. \ref{fig:3}b and c show that the both rectangular tractions within range $ x \in [-x_{0}, x_{0}] $ (where $ x_{0} = 10 {\rm{m}} $ in Table \ref{tab:2}) are zero for the whole time span $ [t_{1}, t_{4}] $, which is identical to the boundary condition in Eq. (\ref{eq:2.21a}). Fig. \ref{fig:3}d shows that no deformation occurs on the ground surface in the range $ x \in [-20,-x_{0}] \cup [x_{0}, 20] $. The deformation of range $ x \in (-\infty, -20) \cup (20,\infty) $ is also zero, but is difficult to illustrate.

In addition, the equality of $ E_{-1} = -R^{2}/2 $ derived by Eqs. (\ref{eq:4.10c}) and (\ref{eq:4.10c'}) is also examined, but the procedure can only be perceived and illustrated by re-running the above-mentioned uploaded codes by the readers themselves. The iterative solution procedure converges very fast, and the maximum iteration rep is $ Q = 42 $ with elapsing time less than 2 secs.

The above results show that the present solution meet the required mixed boundary conditions along ground surface for the whole time span, and further verification can be conducted.

\subsection{Comparisons with finite element solution}%
\label{sub:Comparisons with finite element solution}

\begin{figure}[htpb]
  \begin{center}
    \includegraphics[width = 1.\textwidth]{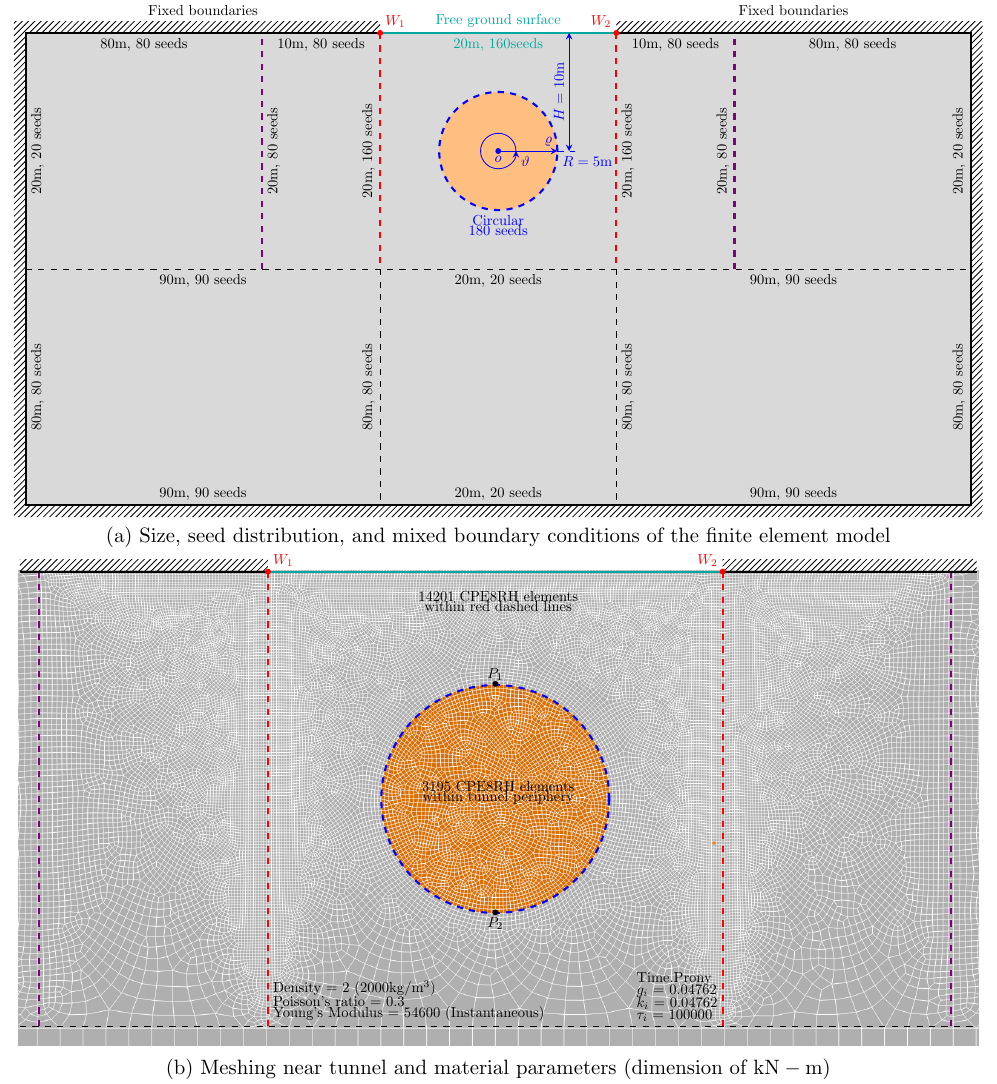}
  \end{center}
  \caption{Finite element model for verification}%
  \label{fig:4}
\end{figure}

The corresponding finite element solution is conducted on \texttt{ABAQUS 2016}, with dimensions of $ {\rm{kN}} $ and $ {\rm{m}} $ for force and length, respectively. The size and seed distributions of the finite element model are illustrated in Fig. \ref{fig:4}a, where the fixed ground surface is much larger than the free one to approximate the infinite fixed ground surface in the analytical solution. Instead of using a corresponding axisymmetrical half, the full finite element model is delibrately used for efficiency comparison between present solution and the finite element one. Meshing near tunnel and the input material parameters are shown in Fig. \ref{fig:4}b. The Young's modulus is computed via $ 2(1+\nu)(G_{\infty}+G_{E}) $ for instantaneous elasticity. The \texttt{Viscoelastic} parameters uses \texttt{Prony} series of \texttt{Time} Domain. The tunnel vault and bottom points respectively marked as $ P_{1} $ and $ P_{2} $ are used to track stress and displacement variation against time span. The steps of finite element solution is illustrated in Table \ref{tab:2}.

\begin{table}[htpb]
  \centering
  \caption{Steps of the finite element solution}
  \label{tab:2}
  \small
  \begin{tabular}{ccccl}
    \toprule
    Step & Procedure & Time period & Incrementation & Loads and Interaction \\
    \midrule
    0 & (Initial) & - & - & Applying fixed boundaries and geostress \\
    1 & Geostatic & 1 & Automatic & Applying gravititational stress field \\
    2 & Static, General & 1 & Automatic & Deactivating element within tunnel periphery \\
    3 & Static, General & 1 & Automatic & Applying static tractions in Eq. (\ref{eq:6.1}) \\
    4 & Visco & 20 & 0.05 & Applying time-dependent tractions in Eq. (\ref{eq:6.2}) \\
    \bottomrule
  \end{tabular}
\end{table}

The boundary condition along tunnel periphery in the finite element solution should be in local polar formation with a local cylindrical coordinate system $ \varrho o \vartheta $ at the tunnel center, as shown in Fig. \ref{fig:4}a. The radial and tangential static tractions along tunnel periphery on the remaining geomaterial side can be expressed as 
\begin{equation}
  \label{eq:6.1}
  \left\{
    \begin{aligned}
      & \sigma_{R}(\vartheta) = \gamma \left( H - R \sin\vartheta \right) \cdot \left( \frac{1+k_{0}}{2} - \frac{1-k_{0}}{2}\cos2\vartheta \right) \\
      & \tau_{R}(\vartheta) = \gamma \left( H - R \sin\vartheta \right) \cdot \frac{1-k_{0}}{2} \sin2\vartheta \\
    \end{aligned}
  \right.
\end{equation}
The two stress components in Eq. (\ref{eq:6.1}) are applied to the \texttt{Analytical Fields} submodule. The time-dependent tractions to simulate progressive shallow tunnelling of tunnel face can be expressed as
\begin{equation}
  \label{eq:6.2}
  \left\{
    \begin{aligned}
      & \sigma_{R}(\vartheta,t) = - U(t) \cdot \sigma_{R}(\vartheta) \\
      & \tau_{R}(\vartheta,t) = - U(t) \cdot \tau_{R}(\vartheta)
    \end{aligned}
  \right.
\end{equation}

The specified equivalent coefficient in Eq. (\ref{eq:2.5'}) is discreted as $ U(t_{i}-t_{1}) $, where $ t_{i} = i\Delta{t} + t_{1} $, and $ \Delta{t} = 0.05 $ for the incrementation of Step 4 in Table \ref{tab:2}. The discreted substeps $i (1 \leq i \leq 400)$ of Step 4 are organized into \texttt{Tabular} formation of $ i\Delta{t} $ and $ U(i\Delta{t}) $ in \texttt{Amplitutdes} submodule.

\begin{figure}[htpb]
  \centering
  \includegraphics[width = 1.\textwidth]{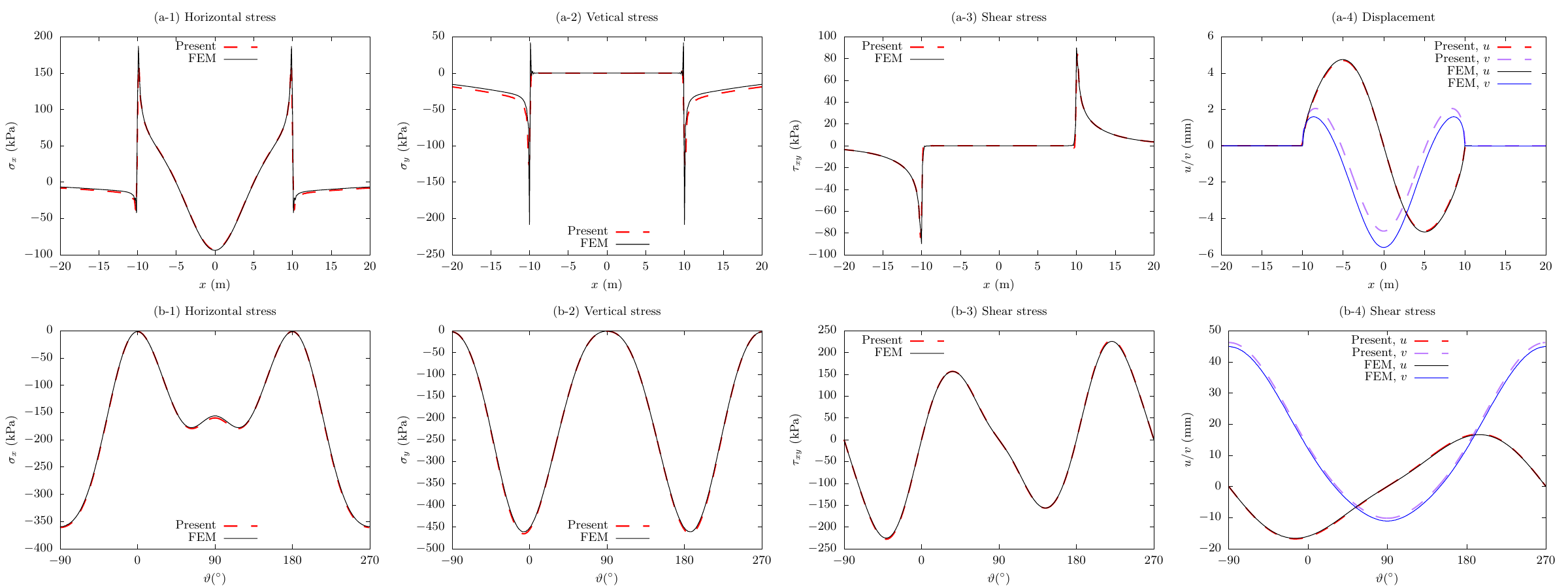}
  \caption{Comparisons of stress and displacement components along ground surface (a) and tunnel periphery (b) for time moment $ t = t_{4} $}
  \label{fig:5}
\end{figure}

With the input parameters above, the finite element results can be obtained. The stress and displacement comparisons along ground surface and tunnel periphery for a given time moment $ t_{4} $ are shown in Figs. \ref{fig:5}a and \ref{fig:5}b, respectively. Both figures illustrate good agreements between the present solution and FEM one, except for a small discrepancy of the vertical displacement along ground surface within ground range $ x \in [-x_{0},x_{0}] $ in Fig. \ref{fig:5}a-4, which is caused by the finite size of $ 200 {\rm{m}} \times 100 {\rm{m}} $ of the finite element model. If the model size tends to infinity, such a displacement discrepancy is expected to be eliminated.

\begin{figure}[htpb]
  \centering
  \includegraphics[width = 1.\textwidth]{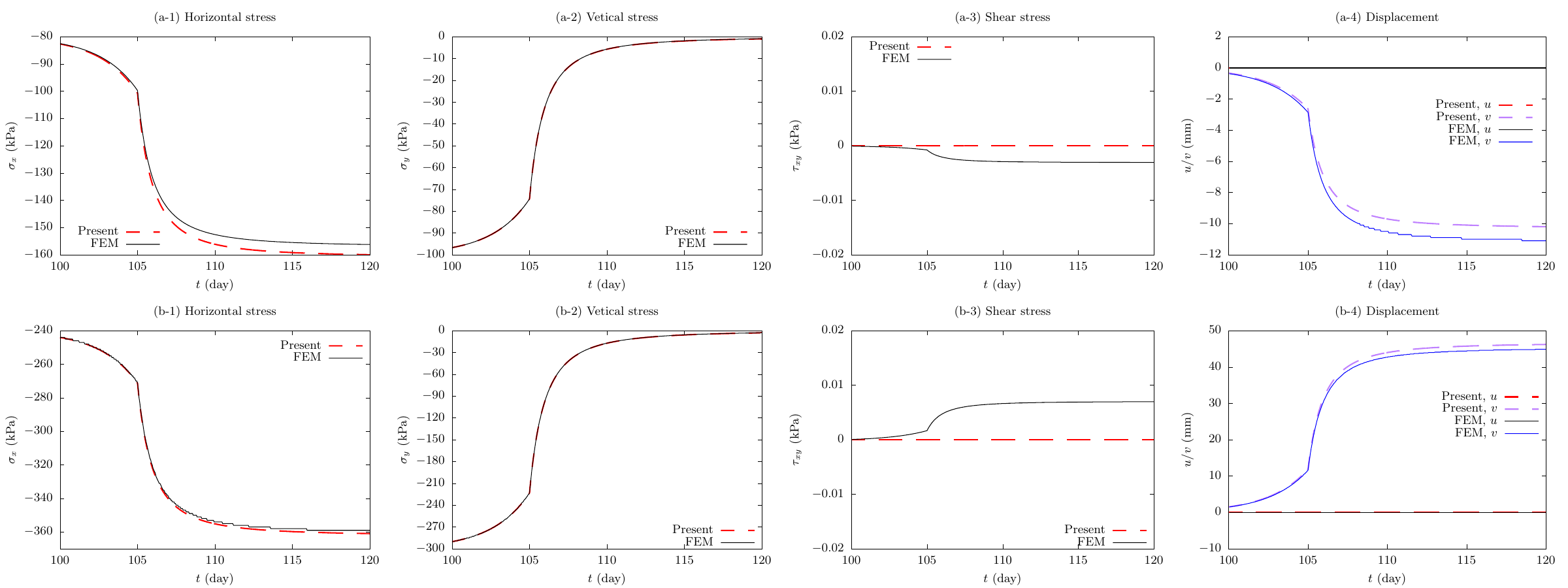}
  \caption{Comparisons of stress and displacement components of points $ P_{1} $ and $ P_{2} $ at tunnel vault for time span $ t \in [ t_{1}, t_{4} ] $}
  \label{fig:6}
\end{figure}

The time-dependent stress and displacement comparisons for points $ P_{1} $ and $ P_{2} $ are shown in Figs. \ref{fig:6}a and \ref{fig:6}b. These two figures suggest good agreements between the present solution and FEM one, especially for the vertical stress $ \sigma_{y} $. The zero $ \tau_{xy} $ for both points indicate that the analytical solution should be more accurate, since the theoretical values of $ \tau_{xy} $ at these two points should always be zero. The time-dependt vertical displacement discrepancies for theset two points are similar to Fig. \ref{fig:5}a-4 for the finite model size.

The verification of the comparisons with the finite element solution in Figs. \ref{fig:5} and \ref{fig:6} further reveal a remarkable feature of the present solution that the elapsing time to obtain the results in Figs. \ref{fig:3}, \ref{fig:5}, and \ref{fig:6} of the present solution is much less than the finite element one. The finite element solution is conducted on a computer of Intel i5-3470 CPU, and takes about 130 minutes for computation alone. Meanwhile, the present solution is conducted on a computer of Intel i3-6100T CPU, and takes about 23 secs for computation alone. These two computers are both outdated with comparable performance, while the computation efficiency of the present solution is far better than the finite element one over 330 times. 

Both verifications above indicate that the present solution is accurate and very efficient, and the present solution can be used for further discussions.

\section{Discussions on present solution}%
\label{sec:further discussion}

In this section, discussions will be made to reveal the some new insights of the verified present solution, comparing to our previous studies with reasonable far-field displacement \cite{LIN2024appl_math_model,LIN2024comput_geotech_1,lin2024over-under-excavation}. For generality, the parameters are normalized as $ G_{E}^{\ast} = G_{E}/G_{\infty} $, $ x_{0}^{\ast} = x_{0}/R $, $ H^{\ast} = H/R $, and $ V^{\ast} = V/R ({\rm{day^{-1}}}) $. All the rectangular stress components and the radial displacement $ u_{r} $ along tunnel periphery in physical plane are normalized by $ \gamma H $ and by $ u_{0} = \frac{\gamma H R}{2G_{\infty}} $, respectively.

\subsection{Excavation rate on the effect of tunnel face}%
\label{sub:Excavation rate on the effect of tunnel face}

Comparing to our previous studies related to shallow tunnelling with reasonable far-field displacement \cite{LIN2024appl_math_model,LIN2024comput_geotech_1,lin2024over-under-excavation}, present solution focuses more on the effect of tunnel face, which would be highly affected by the normalized excavation rate $ V^{\ast} ({\rm{day^{-1}}}) $. To investigate the effect of tunnel excavation rate, the following normalized parameters are selected as $ G_{E}^{\ast} = 0.05 $, $ x_{0}^{\ast} = 2 $, $ H^{\ast} = 2 $, so that the mechanical model is the same as that in Fig. \ref{fig:4}. The normalized excavation rate is selected as $ V^{\ast} = 0.1, 0.2, 0.4, 1, 2 $, which covers a wide range of very slow excavation rate and very fast excavation rate. We select Points $ P_{1} $ at tunnel ceiling and $ P_{2} $ at tunnel bottom (see Fig. \ref{fig:4}b) for time-dependent datum extraction. The lateral coefficient and Poisson's ratio are selected as $ k_{0} = 0.8 $ and $ \nu = 0.3 $, just as those in Table \ref{tab:1}.

\begin{figure}[htp]
  \centering
  \includegraphics[width = 1.\textwidth]{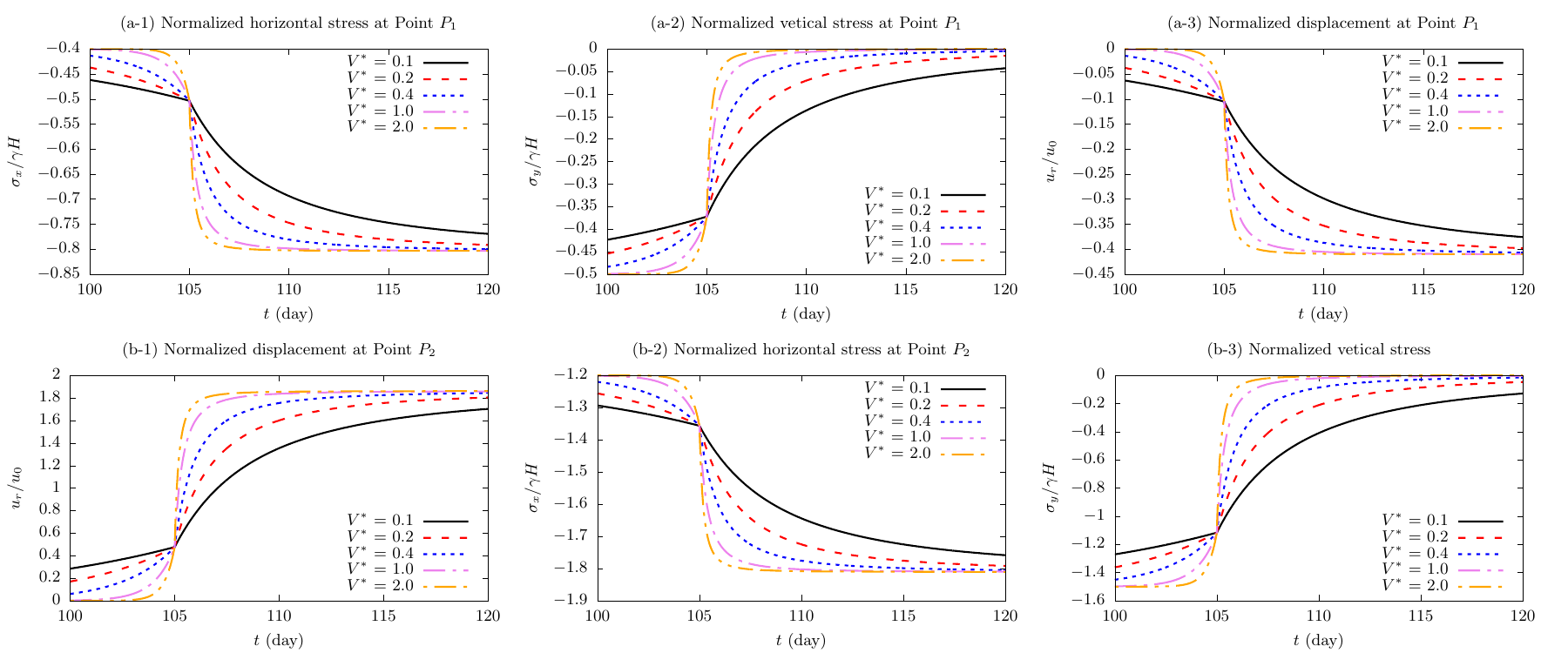}
  \caption{Normalized stress and displacement against excavation time at points $ P_{1} $ and $ P_{2} $ ($ k_{0} = 0.8 $, $ \nu = 0.3 $, $ G_{E}^{\ast} = 0.05 $, $ x_{0}^{\ast} = 2 $, $ H^{\ast} = 2 $)}
  \label{fig:7}
\end{figure}

Subsitituting all necessary parameters into the present solution gives the normalized stress and displacemnt against time for Points $ P_{1} $ and $ P_{2} $ in Fig. \ref{fig:7}. The shear stress and tangential displacement are not presented here, since these two components are zero in theory and in numerical results. Fig. \ref{fig:7} shows that the effect of tunnel face is remarkably reduced when the normalized excavation rate increases from $ V^{\ast} = 0.1 $ to $ V^{\ast} = 2 $, indicating that the effect of tunnel face is not obvious for a fast excavation. The $ V^{\ast} = 0.1 $ curves in Fig. \ref{fig:7} illustrate clear non-smoothness near time moment $ t_{2} = 105 $ (tunnel face), which indicates that the equivalent coefficient $ U(t) $ in Eq. (\ref{eq:2.5}) is probably not suitable for a very slow excavation scheme.

\subsection{Spare shear modulus and viscosity on displacement around tunnel periphery}%
\label{sub:Spare shear modulus}

The spare shear modulus $ G_{E} $ in the Poyting-Thomson model in Eq. (\ref{eq:2.7}) indicates the possible geomaterial weakening due to shallow tunnelling. We should expect that the variation of the normalized spare shear modulus $ G_{E}^{\ast} $ would greatly alter the displacement along tunnel periphery, and we select $ G_{E}^{\ast} = 0, 0.05, 0.1, 0.2, 0.4, 0.5 $ for computation. When $ G_{E}^{\ast} = 0 $, the geomaterial should be degenerated from a viscoelastic one to an elastic one. In our assumption, shallow geomaterial generally shows no rheology, and the Maxwell viscosity parameter takes $ \eta_{E} = 10^{5} {\rm{MPa \cdot day}} $. In this section, we might unlock such a constraint, and take $ \eta_{E} = 10^{2} {\rm{MPa \cdot day}} $ as a possible alternative for theoretical discussion. To be consistent with previous sections, the radial displacements of Points $ P_{1} $ and $ P_{2} $ are selected. The rest parameters take the same values in Table \ref{tab:1}.

\begin{figure}[htp]
  \centering
  \includegraphics[width = 1.\textwidth]{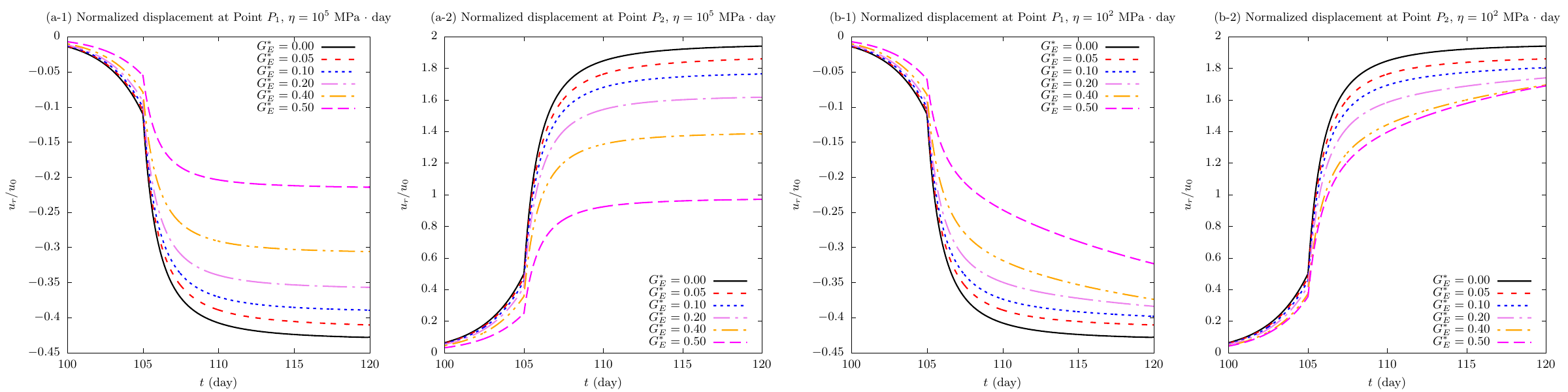}
  \caption{Normalized displacement comparisons at points $ P_{1} $ and $ P_{2} $ ($ k_{0} = 0.8 $, $ \nu = 0.3 $, $ x_{0}^{\ast} = 2 $, $ H^{\ast} = 2 $, $ V^{\ast} = 0.4 $)}
  \label{fig:8}
\end{figure}

Substituting all necessary parameters into the present solution gives the normalized displacement comparisons in Fig. \ref{fig:8}. Fig. \ref{fig:8} shows two typical results that (1) the normalized radial displacements at Points $ P_{1} $ and $ P_{2} $ present overall decreases with a larger $ G_{E}^{\ast} $ due to a larger initial shear modulus $ G_{E} + G_{\infty} $, indicating that the spare shear modulus $ G_{E} $ is more dominant than the viscosity $ \eta_{E} $ in the convolution in Eq. (\ref{eq:5.3}); (2) lower viscosity would cause larger normalized radial displacement, indicating that geomaterial rheology would remarkably increase time-dependent displacement.

\subsection{Solution convergence and paradox}%
\label{sub:Solution convergence and paradox}

\begin{figure}[htp]
  \centering
  \includegraphics[width = 1.\textwidth]{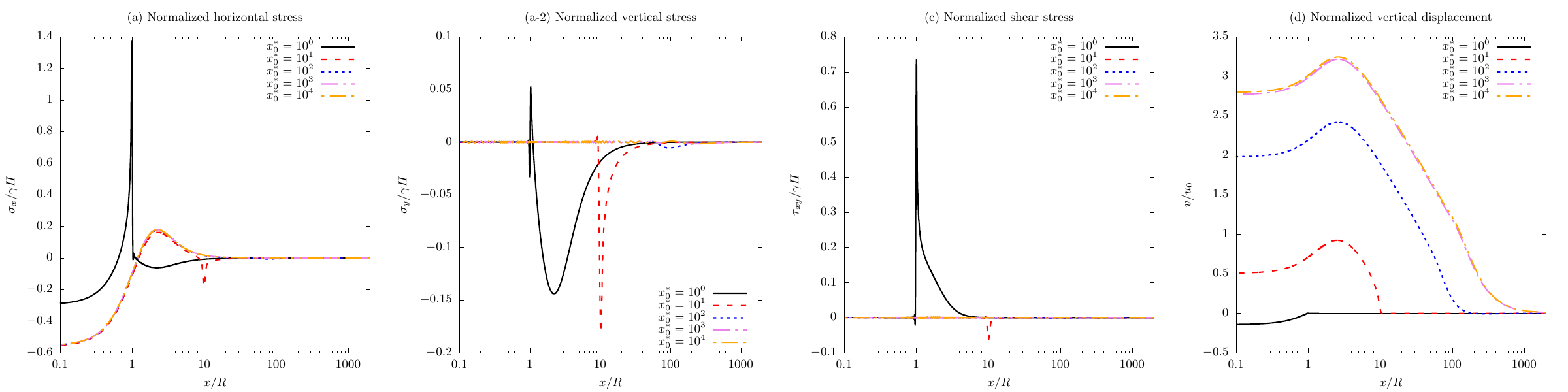}
  \caption{Normalized stress and displacement along ground surface at time moment $ t_{4} $ on variation of fixed ground surface ($ k_{0} = 0.8 $, $ \nu = 0.3 $, $ G_{E}^{\ast} = 0.05 $, $ H^{\ast} = 2 $, $ V^{\ast} = 0.4 $)}
  \label{fig:9}
\end{figure}

In our previous studies \cite{LIN2024appl_math_model,LIN2024comput_geotech_1,lin2024over-under-excavation}, the stress and displacement components along ground surface would be to convergent if the horizontal coordinate $ x_{0} $ of the joint points $ W_{1}(-x_{0},0) $ and $ W_{2}(x_{0},0) $ trends to infinity ($ x_{0} \rightarrow \infty $). In present solution, similar convergence is expected to be observed, especially at time moment $ t_{4} $, where the effects of tunnel face are not felt any more (see Fig. \ref{fig:3}). We select $ x_{0}^{\ast} = 10^{0}, 10^{1}, 10^{2}, 10^{3}, 10^{4} $, and keep the rest parameters the same in Table \ref{tab:1}. The value $ x_{0}^{\ast} = 10^{4} $ is large enough, and can numerically treated as infinity. Substituting all necessary parameters into the present solution gives the stress and displacement components along ground surface of time moment $ t_{4} $ in Fig. \ref{fig:9}. Fig. \ref{fig:9} shows that as the normalized parameter $ x_{0}^{\ast} $ gets larger, all the stress and displacement components trend to clean convergence, which is favorable in mathematics. However, a latent paradox in Fig. \ref{fig:9} should be explained.

When $ x_{0}^{\ast} \rightarrow \infty $, the traction components along ground surface in Figs. \ref{fig:9}b and \ref{fig:9}c meet the possible engineering assumption of free ground surface, while the convergent upheaval displacement in Fig. \ref{fig:9}d is paradoxical to in-situ settlement troughs \cite{Loganathan_Poulos1998,park2004,park2005}, which are overall downward. The reason for such a paradox might be attributed to the difference between the elastic geomaterial assumed in present solution and real-world geomaterials (soils, rocks, and etc.). One typical reason of the upheaval is that the elastic geomaterial is non-porous, and the unbalanced resultant due to excavation is fully applied along tunnel periphery. The other reason is that the periphery of shallow tunelling is assumed to be an ideal circle in present solution, and no over excavation is considered. In contrast, real-world soils are generally porous materials, and over excavation in shallow tunnelling generally exists, both of which cause the unbalanced resultant along a tunnel periphery to be consumed by squeezing the porous voids with certain convergence patterns \cite{Sagaseta_image_method,Verruijt1996polar,park2004,park2005,zhang2017analytical,zhang2018complex}, instead of bringing an upheaval. Though the paradox can be explained with good reasons, we should emphasize that the present solution has certain limits.

\section{Conclusions}%
\label{sec:Conclusions}

A new time-dependent complex variable method on shallow tunnelling in gravitational geomaterial is proposed in this paper, and the far-field displacement singularity due to nonzero resultant along tunnel periphery is eliminated by fixed far-field ground surface. The effect of tunnel face is simplified into a plane strain problem with progressive virtual traction along tunnel periphery, which accompanies the mixed boundary conditions along ground surface to be tranformed into a homogenerous Riemann-Hilbert problem with extra constraints. The problem is subsequently solved in an iterative manner to reach numerically stable stress and displacement solutions. The proposed solution is stepwisely verified via a numerical case by examining the zero values of the mixed boundary conditions, and by comparing with a corresponding finite element solution. Both verifications illustrate good agreements. Further discussions reveal that the effect of tunnel face shows nagative correlation to excavation rate, and lower viscosity would cause larger progressive displacement in geomaterial. In addition, solution convergence paradoxically discloses the limits of our proposed solution due to the non-porosity geomaterial assumption.

\section*{Acknowledgement}

This study is financially supported by the Fujian Provincial Natural Science Foundation of China (Grant No. 2022J05190 and 2023J01938), the Scientific Research Foundation of Fujian University of Technology (Grant No. GY-Z21026), and the National Natural Science Foundation of China (Grant No. 52178318).

%% \section{}
%% \label{}

%% If you have bibdatabase file and want bibtex to generate the
%% bibitems, please use
%%
%%  \bibliographystyle{elsarticle-num} 
%%  \bibliography{<your bibdatabase>}

%% else use the following coding to input the bibitems directly in the
%% TeX file.
\bibliographystyle{elsarticle-num}
\bibliography{newbib}

% \begin{thebibliography}{00}

%% \bibitem{label}
%% Text of bibliographic item

% \bibitem{}

% \end{thebibliography}

\clearpage

\clearpage

\clearpage

\clearpage

\clearpage

\clearpage

\clearpage

\clearpage

\clearpage

\end{document}